\documentclass{amsart}
\usepackage{fullpage,url}
\usepackage{graphicx}
\usepackage{amsmath}
\usepackage{amssymb}
\usepackage{amsfonts}
\usepackage{amsthm}
\usepackage{enumerate}
\usepackage{tkz-graph,tkz-berge}

\newcommand{\Z}{\mathbb{Z}}
\renewcommand{\P}{\mathcal{P}}

\newcommand{\ra}{\rightarrow}

\newcommand{\<}{\langle}
\renewcommand{\>}{\rangle}
\renewcommand{\l}{\ell}
\newcommand{\lt}{\leadsto}

\begin{document}

\newtheorem{Theo}{Theorem}[section]
\newtheorem{Prop}[Theo]{Proposition}
\newtheorem{Lemma}[Theo]{Lemma}
\newtheorem{Cor}[Theo]{Corollary}

\theoremstyle{definition}
\newtheorem{Defn}[Theo]{Definition}
\newtheorem{Remark}[Theo]{Remark}
\newtheorem{Exam}[Theo]{Example}

\title{Graph Invertibility}
\author{Cam McLeman and Erin McNicholas}
\begin{abstract}
Extending the work of Godsil and others, we investigate the notion of
the inverse of a graph (specifically, of bipartite graphs with a
unique perfect matching).  We provide a concise necessary and
sufficient condition for the invertibility of such graphs and
generalize the notion of invertibility to multigraphs. We examine the
question of whether there exists a ``litmus subgraph'' whose
bipartiteness determines invertibility. As an application of our
invertibility criteria, we quickly describe all invertible unicyclic
graphs. Finally, we describe a general combinatorial procedure for
iteratively constructing invertible graphs, giving rise to large new
families of such graphs.
\end{abstract}

\maketitle

\section{Intro}
Given the plethora of composition operations on
graphs\footnote{``Graphs'' in this article can include multiple edges
  between distinct vertices, but no loops.  We will occasionally use
  the term multigraph when directly contrasting results to the
  corresponding properties for simple graphs.  } (Cartesian sum,
tensor product, etc.), one is naturally led to the question of whether
or not there is a sensible notion of the \emph{inverse} of a graph.
There is no shortage of possible definitions: A first attempt is to
define two graphs to be inverses if they possess inverse adjacency
matrices.  This turns out to be overly restrictive, as under this
definition only the graphs $nK_2$ are invertible, with themselves as
their own inverses (\cite{HM76}).  A second attempt, motivated by the
observation that the eigenvalues of the sum and product of two graphs
are the pairwise sums and products of the eigenvalues of the original
graphs, is to call a graph $G$ invertible if there exists another
graph $G^{-1}$ such that for each eigenvalue $\lambda$ of $G$,
$\frac{1}{\lambda}$ is an eigenvalue of $G^{-1}$ (with the same
multiplicity).  This definition too allows some unfortunate phenomena:
If $G_1$ and $G_2$ are cospectral and non-isomorphic, then $G_1^{-1}$
and $G_2^{-1}$ (if such graphs exist) both satisfy the criterion for
being inverses to $G_1$, and we are left with multiple non-isomorphic
inverses.  Further, there would be no hope of attaining the obviously
desirable property that $(G^{-1})^{-1}$ be isomorphic to $G$.

\bigskip

It therefore behooves us to strengthen the condition defining the
inverse.  We begin by noting that since adjacency matrices are
diagonalizable (being real and symmetric), two such matrices are
cospectral if and only if they are similar.  The reciprocal eigenvalue
condition described above is thus tantamount to asserting that the
inverse $A^{-1}$ of the adjacency matrix to $G$ is similar to the
adjacency matrix of $G^{-1}$.  A strengthening of the definition comes
from a result of Godsil (\cite{Go85}) that under certain conditions on
$G$ (described below), the inverse adjacency matrix $A^{-1}$ is in
fact \emph{signable} to a non-negative symmetric integral matrix with
zeros on the diagonal, i.e., to the adjacency matrix of a graph.  Here
we say $A$ is signable to $B$ if $A$ can be conjugated to $B$ by a
diagonal matrix whose diagonal entries are all $\pm 1$ (i.e., by a
\emph{signing matrix}).  We therefore adopt the following definition:

\begin{Defn}
Given a graph $G$, we say that a graph $H$ is an \emph{inverse} of $G$
if they possess adjacency matrices $A_G$ and $A_H$ such that $A_H$ is
signable to $A_G^{-1}$.  We then say that $G$ is \emph{invertible},
and say that $G$ is \emph{simply invertible} if there exists a simple
graph $H$ which is an inverse of $G$.  (In particular, we emphasize
that a simple graph can be invertible but not simply invertible.)
\end{Defn}

Clearly this stronger condition defining invertibility implies the
earlier reciprocal eigenvalue property, and it is thus easy to find
non-invertible graphs -- namely, any graph with an eigenvalue of 0,
e.g., bipartite graphs on an odd number of vertices.  In fact, this is
a convenient place to note that for an invertible graph $G$ with an
inverse $H$, we must have
$\det(A_G)^{-1}=\det(A_G^{-1})=\det(A_H)\in\Z$, and so $\det(A_G)=\pm
1$ for any invertible graph.  This forces $G$ to admit a perfect
matching (or ``1-factor''), providing fairly compelling evidence that
most graphs are not invertible.  Following Godsil and the
subsequent literature, we focus on graphs $G$ which are bipartite and
have a \emph{unique} perfect matching $M$.  The first significant
invertibility result (\cite{Go85}, Theorem 2.2) gives that a simple
graph $G$ (bipartite with a unique perfect matching $M$) is invertible
if the graph $G/M$ obtained by contracting each edge of $M$ is
bipartite.  The aim of the current paper is to extend results of this
form in a variety of different directions.

\subsubsection*{Summary of Results}

Section 2 contains preliminaries on bipartite graphs with a unique
perfect matching, focusing on inversion and extending previously
well-known results for simple graphs to the context of multigraphs.
In particular, we give a purely graph-theoretic construction of the
inverse (when it exists -- see Theorem \ref{gottabe}), which we dub
the \emph{parity closure} of the graph.  We emphasize a graphical
point of view (as opposed to a poset-theoretical or linear-algebraic
one), enough so that it is frequently possible to bypass any
matrix-inversion calculations and ``eyeball'' both the invertibility
and inverse of a given graph.  Further, we prove that the construction
satisfies the desired properties of an inverse from the introduction
(i.e., that $(G^{-1})^{-1}=G$ -- see Theorem \ref{G++}).  In Section
3, we turn our attention to determining conditions for the inverse to
exist.  First, we extend a variety of known results on invertibility
to the context of multigraphs, among them the result of Godsil
mentioned above and a related result of \cite{TK09} that a necessary
condition for invertibility is the bipartiteness of a certain subgraph
$\Gamma$ of $G/M$.  Continuing, we note that the main result of
\cite{TK09} gives much more, reducing the question of the
invertibility of $G$ to the invertibility of a collection of subgraphs
(the ``undirected intervals'') of $G$.  Their culminating necessary
and sufficient condition for invertibility admits some curiosities,
however. If $G$ is either:

\smallskip

\begin{enumerate}[(a)]
\item A simply invertible undirected interval graph with bipartite
  Hasse diagram (Figure \ref{TKfigs}, left); or
\item A non-invertible interval graph with bipartite Hasse diagram all
  of whose proper sub-undirected intervals are invertible (Figure \ref{TKfigs},
  right),
\end{enumerate}

\begin{figure}[!ht]
\begin{center}
\begin{tabular}{ccc}
\scalebox{.5}{
\begin{tikzpicture}
\GraphInit[vstyle=Normal]
\SetVertexMath
\SetVertexNoLabel

\Vertex[x=0.0cm,y=0.0cm]{a0}
\Vertex[x=1.6667cm,y=0.0cm]{a1}
\Vertex[x=3.3333cm,y=0.0cm]{a2}
\Vertex[x=5.0cm,y=0.0cm]{a3}
\Vertex[x=0.0cm,y=2.5cm]{a4}
\Vertex[x=1.6667cm,y=2.5cm]{a5}
\Vertex[x=3.3333cm,y=2.5cm]{a6}
\Vertex[x=5.0cm,y=2.5cm]{a7}
\Vertex[x=2.5cm,y=0.75cm]{a8}
\Vertex[x=2.5cm,y=1.75cm]{a9}
\Vertex[x=1.6667cm,y=3.75cm]{a10}
\Vertex[x=3.3333cm,y=3.75cm]{a11}

\Edge(a0)(a1)
\Edge(a2)(a3)
\Edge(a2)(a6)
\Edge(a2)(a8)
\Edge(a4)(a5)
\Edge(a5)(a6)
\Edge(a6)(a7)
\Edge(a6)(a9)
\renewcommand*{\EdgeLineWidth}{1mm}
\tikzset{EdgeStyle/.style   = {line width = \EdgeLineWidth}}
\Edge(a3)(a7)
\Edge(a8)(a9)
\Edge(a0)(a4)
\Edge(a1)(a2)
\Edge(a5)(a10)
\Edge(a6)(a11)

\end{tikzpicture}}
&
\hspace*{.5in}
&
\scalebox{.5}{
\begin{tikzpicture}
\GraphInit[vstyle=Normal]
\SetVertexMath
\SetVertexNoLabel

\Vertex[x=0.0cm,y=0.0cm]{a0}
\Vertex[x=1.6667cm,y=0.0cm]{a1}
\Vertex[x=3.3333cm,y=0.0cm]{a2}
\Vertex[x=5.0cm,y=0.0cm]{a3}
\Vertex[x=0.0cm,y=2.5cm]{a4}
\Vertex[x=1.6667cm,y=2.5cm]{a5}
\Vertex[x=3.3333cm,y=2.5cm]{a6}
\Vertex[x=5.0cm,y=2.5cm]{a7}
\Vertex[x=1.6667cm,y=3.75cm]{a10}
\Vertex[x=3.3333cm,y=3.75cm]{a11}

\Edge(a0)(a1)
\Edge(a2)(a3)
\Edge(a2)(a6)
\Edge(a4)(a5)
\Edge(a5)(a6)
\Edge(a6)(a7)
\renewcommand*{\EdgeLineWidth}{1mm}
\tikzset{EdgeStyle/.style   = {line width = \EdgeLineWidth}}
\Edge(a3)(a7)
\Edge(a0)(a4)
\Edge(a1)(a2)
\Edge(a5)(a10)
\Edge(a6)(a11)

\end{tikzpicture}}
\end{tabular}
\end{center}
\caption{}\label{TKfigs}
\end{figure}
\noindent then the principal result (\cite{TK09}, Theorem 2.6) returns
a tautology -- $G$ is invertible if and only if $G$ is invertible.
The main result of Section 3 (Theorem \ref{Alg}) replaces undirected
intervals with different key substructures which completely determine
the invertibility of the graph, leading to a concise necessary and
sufficient condition.  In particular, Examples \ref{introleft} and
\ref{introright} determine the invertibility of the graphs in Figure \ref{TKfigs}
via a trivial calculation (especially in comparison to inverting and
correctly signing a $12\times 12$ or $10\times 10$ matrix).  As a more
substantial application, we quickly recover the characterization of
invertible unicyclic graphs found in \cite{AK07}.  We close the
section focusing on the striking appearance of bipartiteness in the
inversion results of both \cite{Go85} and \cite{TK09} and investigate
the question of whether there exists an ``optimal'' subgraph of $G/M$
in the sense that its bipartiteness is equivalent to the invertibility
of $G$.  We answer this in the negative but improve both previous
results in the sense that we find a subgraph of $G/M$ whose
bipartiteness implies the invertibility of $G$, and a supergraph of
$\Gamma$ whose bipartiteness is implied by the invertibility of $G$
(Theorems \ref{minimalprimepairs} and \ref{DeltaD}, respectively).

\medskip

Finally, Section 4 addresses the question of how to construct
invertible bipartite graphs with a unique perfect matching.  If one
views a graph as being constructed via a sequence of operations
consisting of adding a new vertex and edges to that vertex, then it is
natural to ask when such an operation preserves the invertibility of
the graph.  A complete answer to this question would provide a purely
combinatorial description of the class of invertible graphs.  Theorem
\ref{IneqThm} gives a necessary and sufficient condition for such an
operation to preserve the invertibility of the graph, providing a
rather large array of constructible classes of invertible graphs (see,
e.g., Proposition \ref{2setprop}.)  As a demonstration of the utility,
we return in Section 4.1 to the topic of unicyclic graphs, and
describe explicit combinatorial constructions of invertible unicyclic
graphs with prescribed size and cycle length.

%*****************************************************************

\section{Preliminaries}\label{prelims}

%*****************************************************************

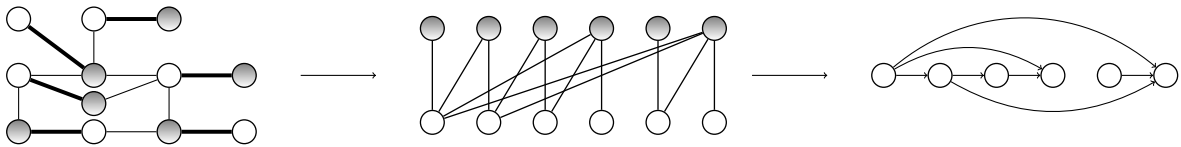
\begin{figure}[!ht]
\begin{center}
\scalebox{.5}{
\begin{tikzpicture}
\GraphInit[vstyle=Normal]
\SetVertexMath
\SetVertexNoLabel

%Top Left

\Vertex[x=-3cm,y=-.25cm,style={shade}]{a8}
\Vertex[x=-1cm,y=-.25cm]{a2}
\Vertex[x=1cm,y=-.25cm,style={shade}]{a9}
\Vertex[x=3cm,y=-.25cm]{a3}
\Vertex[x=-3cm,y=1.25cm]{a1}
\Vertex[x=-1cm,y=1.25cm,style={shade}]{a11}
\Vertex[x=1cm,y=1.25cm]{a0}
\Vertex[x=3cm,y=1.25cm,style={shade}]{a6}
\Vertex[x=-3cm,y=2.75cm]{a5}
\Vertex[x=-1cm,y=2.75cm]{a4}
\Vertex[x=1cm,y=2.75cm,style={shade}]{a10}
\Vertex[x=-1cm,y=.5cm,style={shade}]{a7}

\Edge(a4)(a11)
\Edge(a1)(a11)
\Edge(a11)(a0)
\Edge(a0)(a9)
\Edge(a1)(a8)
\Edge(a2)(a9)
\Edge(a7)(a0)

\tikzset{EdgeStyle/.style   = {line width = 3pt}}
\Edge(a0)(a6)
\Edge(a1)(a7)
\Edge(a2)(a8)
\Edge(a9)(a3)
\Edge(a5)(a11)
\Edge(a4)(a10)

%Top Middle
\Vertex[x=8.0cm,y=0.0cm]{b0}
\Vertex[x=9.5cm,y=0.0cm]{b1}
\Vertex[x=11cm,y=0.0cm]{b2}
\Vertex[x=12.5cm,y=0.0cm]{b3}
\Vertex[x=14cm,y=0cm]{b4}
\Vertex[x=15.5cm,y=0cm]{b5}
\Vertex[x=8.0cm,y=2.5cm,style={shade}]{b6}
\Vertex[x=9.5cm,y=2.5cm,style={shade}]{b7}
\Vertex[x=11cm,y=2.5cm,style={shade}]{b8}
\Vertex[x=12.5cm,y=2.5cm,style={shade}]{b9}
\Vertex[x=14cm,y=2.5cm,style={shade}]{b10}
\Vertex[x=15.5cm,y=2.5cm,style={shade}]{b11}
\tikzset{EdgeStyle/.style   = {line width = 1pt}}
\Edge(b4)(b11)
\Edge(b1)(b11)
\Edge(b11)(b0)
\Edge(b0)(b9)
\Edge(b1)(b8)
\Edge(b2)(b9)
\Edge(b7)(b0)

\Edge(b0)(b6)
\Edge(b1)(b7)
\Edge(b2)(b8)
\Edge(b3)(b9)
\Edge(b4)(b10)
\Edge(b5)(b11)

%Top Right
\Vertex[x=20.0cm,y=1.25cm]{c0}
\Vertex[x=21.5cm,y=1.25cm]{c1}
\Vertex[x=23cm,y=1.25cm]{c2}
\Vertex[x=24.5cm,y=1.25cm]{c3}
\Vertex[x=26cm,y=1.25cm]{c4}
\Vertex[x=27.5cm,y=1.25cm]{c5}

\tikzset{EdgeStyle/.style = {thick}}
\tikzset{EdgeStyle/.append style = {->}}

\Edge(c0)(c1)
\Edge(c1)(c2)
\Edge(c2)(c3)
\Edge(c4)(c5)

\tikzset{EdgeStyle/.style = {bend right}}
\tikzset{EdgeStyle/.append style = {->}}
\Edge(c1)(c5)

\tikzset{EdgeStyle/.style = {bend left}}
\tikzset{EdgeStyle/.append style = {->}}
\Edge(c0)(c3)
\Edge[style={bend left=40}](c0)(c5)

%Arrows
\draw[->,thick](4.5,1.25) -- ( 6.5,1.25);
\draw[->,thick](16.5,1.25) -- (18.5,1.25);
\end{tikzpicture}}
\end{center}
\caption{Construction of the associated digraph.}\label{digraphpicture}
\end{figure}

As mentioned in the introduction, we restrict our attention to
bipartite (multi-)graphs on $2n$ vertices with a unique perfect
matching.  We begin with a reduction process to simplify the
discussion.  First, any such graph can be visualized by arranging the
edges of the perfect matching as vertical columns, with orientation
chosen so that the top of each column is the same color (under some
proper 2-coloring with colors, say, black and white).  Next, uniformly
orienting diagonal edges from black vertices to white vertices, we
choose a topological sort of the columns so that the diagonal edges
all have ``positive slope,'' as in the figure below.  Finally, we
construct the \emph{digraph $D=D_G$ associated to $G$} by collapsing
each vertical column to a single vertex.  More precisely, $D_G$ is the
digraph whose vertices are the edges of the perfect matching with an
edge from the $i$-th vertex to the $j$-th vertex of $D_G$ if there is
an edge from the bottom of the $i$-th column to the top of the $j$-th
column in $G$.  Figure \ref{digraphpicture} illustrates the
construction of the associated digraph of a bipartite graph whose
unique perfect matching is drawn in bold.

\medskip

Note that the choice of topological sort will fix once and for all an
ordering of the vertices of $D$, which we will consistently label by
$\{1,2,\ldots,n\}$.  If we now impose the labeling on $G$ where the
bottom row of vertices of $G$ is $\{1,2,\ldots,n\}$ and the top row of
vertices is $\{n+1,\ldots,2n\}$, then the adjacency matrix $A=A_G$ of
$G$ has the simple block form

$$A=\left[\begin{array}{cc}0 & B \\ B^T & 0\end{array}\right],$$
\noindent where $B$ is the $n\times n$ matrix given by $B_{ii}=1$ for
all $i$, and for $i\neq j$, $B_{ij}$ is the number of edges between
columns $i$ and $j$ in $G$.  The point of this construction is that
the associated digraph $D$ is a simpler object than $G$, yet contains
all the information germane to its invertibility.  Namely, since we
have removed the matching edges when constructing $D$, its adjacency
matrix is given by the upper-triangular matrix $B-I$.  We have the
following:

\begin{Theo}\label{intInv}
Given a bipartite multigraph $G$ with unique perfect matching $M$, its
adjacency matrix $A$ is invertible with integral inverse
matrix

\begin{align*}
A^{-1}=\left[\begin{array}{cc}0 & (B^T)^{-1} \\ B^{-1} & 0\end{array}\right].
\end{align*}
\end{Theo}
\begin{proof}
The form of the inverse is readily verified by matrix multiplication.
For the integrality condition, note that $B$ is upper triangular with
ones on the diagonal and so has determinant 1.  By the inverse-adjugate
formula, all the entries of $B^{-1}$ are thus integers.
\end{proof}
Recall from the introduction that by definition $G$ is invertible if
and only if $A^{-1}$ is signable to the adjacency matrix of another
graph, i.e., if there exists a (diagonal) signing matrix $S=(s_{ii})$,
with $s_{ii}=\pm 1$ for all $1\leq i\leq n$, such that $SA^{-1}S$ is
an integral symmetric matrix with non-negative entries.  Of course,
since $SA^{-1}S$ is necessarily symmetric and integral, and
$|(SA^{-1}S)_{ij}|=|(A^{-1})_{ij}|$, it suffices to check the existence
of a signing matrix such that $SA^{-1}S=|A^{-1}|$, where $|A^{-1}|$
denotes the matrix obtained by taking absolute values of $A^{-1}$
componentwise.  Since we can conjugate block matrices by blocks, the
matrix $A^{-1}$ is signable if and only the matrix $B^{-1}$ is.  If we
abuse terminology slightly and say that $D$ itself is
\emph{invertible} if $B^{-1}$ is signable, we have now proved the
following:
\begin{Theo}\label{signable}
A bipartite graph $G$ with unique perfect matching is invertible if
and only if its associated digraph $D$ is, i.e., if there exists a
signing matrix $S$ such that $SB^{-1}S=|B^{-1}|$.
\end{Theo}
The problem of inverting $G$ is thus reduced to the problem of signing
$B$, and so we investigate the question of $G$'s invertibility through
the lens of combinatorics on its associated digraph.  As a first
step in this direction, using the nilpotency of $B-I$, we can
calculate the entries of $B^{-1}$ via
\begin{equation*}
B^{-1}=(I+(B-I))^{-1}=\sum_{m=0}^\infty(-1)^m(B-I)^m.
\end{equation*}

By standard results on adjacency matrices, $(B-I)^m_{i,j}$ counts the
number of directed paths of length $m$ between vertices $i$ and $j$ in
$D$, and so the entries of $B^{-1}$ are given by summing this quantity
over all $m\geq 0$, each term weighted by $\pm 1$ according to the
parity of $m$:

\begin{equation}\label{Inv}
B_{i,j}^{-1}=\#(\text{even-length paths from $i$ to $j$ in $D$})-\#(\text{odd-length
  paths from $i$ to $j$ in $D$}).
\end{equation}

Let $\P_{i,j}=\P_{i,j}(D)$ denote the set of all paths between vertices $i$ and
$j$ in a digraph $D$, and $l(P)$ be the length of a path $P$.  Equation
\eqref{Inv} can be rewritten as
\begin{equation}\label{Bij}
B^{-1}_{i,j}=\sum_{P\in\P_{i,j}(D)}(-1)^{l(P)}.
\end{equation}

\begin{Remark}
In the case that $G$ is a \emph{simple} graph, the above are special
cases of poset-theoretical M\"{o}bius-inversion arguments.  Namely,
note that $D$ can be thought of as defining a partial order on the
set of vertices of $G$.  Let $P$ be this poset, and let $\zeta$ and
$\mu$ be the Zeta and M\"{o}bius functions of this poset \cite{Aig79}.
We have $B_{i,j}=\zeta(i,j)$ and $B^{-1}_{i,j}=\mu(i,j)$, and so
inverting $B$ amounts to being able to calculate $\mu(i,j)$ for each
$i\leq j$ between $1$ and $n$.  M\"{o}bius inversion of $\zeta$ then
gives the above formula.
\end{Remark}

Equation \eqref{Bij} hints to an explicit construction of a potential
inverse.  As a motivating example, consider the case that $D$ is a
directed tree, so that there is at most one path between any two
vertices.  In this case, $B_{i,j}^{-1}\in\{0,\pm1 \}$ for all $i$ and
$j$, with a non-zero value if and only if there is a directed path
from $i$ to $j$.  Since any inverse to $D$ should have exactly
$|B_{i,j}^{-1}|$ edges from $i$ to $j$, the upshot of this discussion
is that the inverse should be a graph with an edge from $i$ to $j$ if
and only if $B_{i,j}\neq 0$, i.e., if and only if there is a path from
$i$ to $j$ in $D$.  Such a graph is trivial to construct: In the
figure below, a directed tree $D$ is given with five (solid) edges,
and its inverse is constructed by adding in the dashed edges.

\begin{figure}[!ht]\label{fullexample}
\begin{center}
\scalebox{.6}{
\begin{tikzpicture}
\GraphInit[vstyle=Normal]
\SetVertexMath
\SetVertexNoLabel

\Vertex[x=0.0cm,y=0.0cm]{a0}
\Vertex[x=2.0cm,y=0.0cm]{a1}
\Vertex[x=4.0cm,y=0.0cm]{a2}
\Vertex[x=6.0cm,y=0.0cm]{a3}
\Vertex[x=8.0cm,y=0.0cm]{a4}
\Vertex[x=10.0cm,y=0.0cm]{a5}

\AssignVertexLabel{a}{10}{$1$,$2$,$3$,$4$,$5$,$6$}
\tikzset{EdgeStyle/.append style = {->}}
\Edge(a0)(a1)
\Edge(a2)(a3)
\Edge(a3)(a4)
\Edge(a4)(a5)
\tikzset{EdgeStyle/.style = {->,bend left}}
\Edge(a1)(a5)
\tikzset{EdgeStyle/.style = {->,dashed,bend left}}
\Edge(a2)(a4)
\Edge[style={bend left=40}](a0)(a5)
\tikzset{EdgeStyle/.style = {->,dashed,bend right}}
\Edge(a3)(a5)
\Edge[style={bend left=40}](a2)(a5)
\end{tikzpicture}}
\end{center}
\caption{The inverse of a directed tree.}
\end{figure}
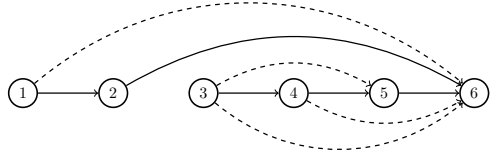

Note that this is the digraph corresponding to the transitive closure
of the original digraph.  Now if $D$ admits (undirected) cycles, then
Equation (3) is more complicated (e.g., two paths from $i$ to $j$ of
opposite parity cancel each other out), but still describes a putative
inverse via a ``parity-corrected'' transitive closure.  We formalize
this in a definition.

\begin{Defn}\label{parclose}
The \emph{parity closure} of a digraph $D$ is the directed graph $D^+$ on
the vertices of $D$ with exactly
\begin{align*}
|B_{i,j}^{-1}|&=|\#(\text{even-length paths from $i$ to $j$ in $D$})-\#(\text{odd-length paths from $i$ to $j$ in $D$})|\\
&=\big|\hspace*{-.15in}\sum_{P\in\P_{i,j}(D)}(-1)^{l(P)}\big|
\end{align*}
edges from $i$ to $j$ for each $i<j$.  If $G$ is a bipartite graph
with unique perfect matching and associated digraph $D$, we define
$G^+$ to be the graph with adjacency matrix
\begin{align*}
A^+=\left[\begin{array}{cc}0 & |(B^T)^{-1}| \\ |B^{-1}| & 0\end{array}\right].
\end{align*}
(Here again, $|B|$ denotes taking the componentwise absolute value of the matrix.)
\end{Defn}

\noindent We summarize the above discussion as the following theorem:
\begin{Theo}\label{gottabe}
If $D$ (resp. $G$) is invertible, then $D^+$ (resp. $G^+$) is its inverse.
\end{Theo}
\begin{Remark}
We remark that if one were to adopt the notation of $G^{-1}$ for the
inverse of $G$ as in the introduction, the previous theorem could be
re-phrased as the equality $G^+=G^{-1}$ for invertible graphs $G$.  In
the sequel we will focus primarily on $G^+$ instead of $G^{-1}$ as the
former is defined for all graphs of interest.
\end{Remark}
\begin{Exam}
We complete in Figure \ref{finishedpicture} the process begun in Figure
\ref{digraphpicture}, constructing the parity closure of a bipartite
graph $G$ with a unique perfect matching.  We form the associated
digraph, take its parity closure (labels indicating multiple edges),
and then return the graph to its original configuration.  The dashed
arrow on the left represents the composite process.

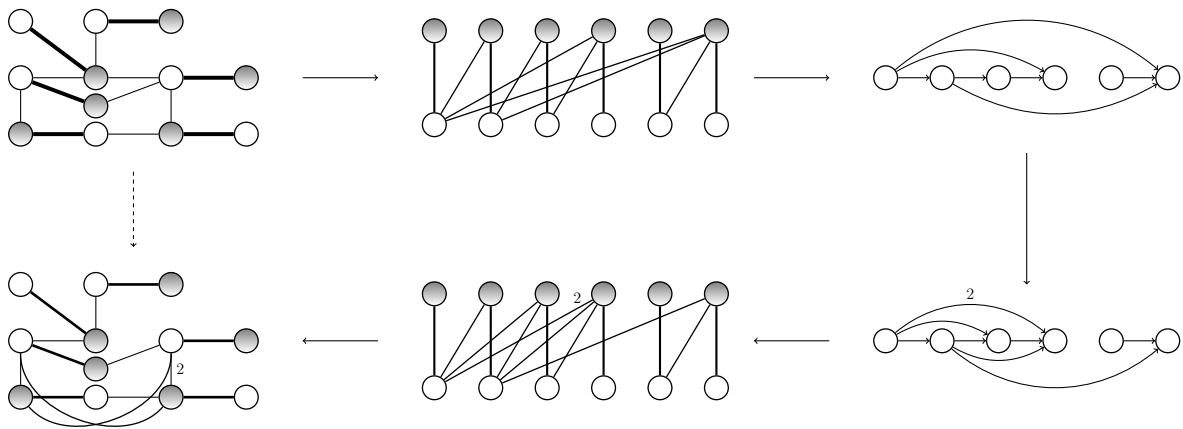
\begin{figure}[!ht]
\begin{center}
\scalebox{.5}{
\begin{tikzpicture}
\GraphInit[vstyle=Normal]
\SetVertexMath
\SetVertexNoLabel

%Top Left
\Vertex[x=-3cm,y=-.25cm,style={shade}]{a8}
\Vertex[x=-1cm,y=-.25cm]{a2}
\Vertex[x=1cm,y=-.25cm,style={shade}]{a9}
\Vertex[x=3cm,y=-.25cm]{a3}
\Vertex[x=-3cm,y=1.25cm]{a1}
\Vertex[x=-1cm,y=1.25cm,style={shade}]{a11}
\Vertex[x=1cm,y=1.25cm]{a0}
\Vertex[x=3cm,y=1.25cm,style={shade}]{a6}
\Vertex[x=-3cm,y=2.75cm]{a5}
\Vertex[x=-1cm,y=2.75cm]{a4}
\Vertex[x=1cm,y=2.75cm,style={shade}]{a10}
\Vertex[x=-1cm,y=.5cm,style={shade}]{a7}

\Edge(a4)(a11)
\Edge(a1)(a11)
\Edge(a11)(a0)
\Edge(a0)(a9)
\Edge(a1)(a8)
\Edge(a2)(a9)
\Edge(a7)(a0)

\tikzset{EdgeStyle/.style   = {line width = 3pt}}
\Edge(a0)(a6)
\Edge(a1)(a7)
\Edge(a2)(a8)
\Edge(a9)(a3)
\Edge(a5)(a11)
\Edge(a4)(a10)
\tikzset{EdgeStyle/.style   = {line width = 1pt}}

%Top Middle
\Vertex[x=8.0cm,y=0.0cm]{b0}
\Vertex[x=9.5cm,y=0.0cm]{b1}
\Vertex[x=11cm,y=0.0cm]{b2}
\Vertex[x=12.5cm,y=0.0cm]{b3}
\Vertex[x=14cm,y=0cm]{b4}
\Vertex[x=15.5cm,y=0cm]{b5}
\Vertex[x=8.0cm,y=2.5cm,style={shade}]{b6}
\Vertex[x=9.5cm,y=2.5cm,style={shade}]{b7}
\Vertex[x=11cm,y=2.5cm,style={shade}]{b8}
\Vertex[x=12.5cm,y=2.5cm,style={shade}]{b9}
\Vertex[x=14cm,y=2.5cm,style={shade}]{b10}
\Vertex[x=15.5cm,y=2.5cm,style={shade}]{b11}

\Edge(b4)(b11)
\Edge(b1)(b11)
\Edge(b11)(b0)
\Edge(b0)(b9)
\Edge(b1)(b8)
\Edge(b2)(b9)
\Edge(b7)(b0)

\tikzset{EdgeStyle/.append style = {ultra thick}}
\Edge(b0)(b6)
\Edge(b1)(b7)
\Edge(b2)(b8)
\Edge(b3)(b9)
\Edge(b4)(b10)
\Edge(b5)(b11)

%Top Right
\Vertex[x=20.0cm,y=1.25cm]{c0}
\Vertex[x=21.5cm,y=1.25cm]{c1}
\Vertex[x=23cm,y=1.25cm]{c2}
\Vertex[x=24.5cm,y=1.25cm]{c3}
\Vertex[x=26cm,y=1.25cm]{c4}
\Vertex[x=27.5cm,y=1.25cm]{c5}

\tikzset{EdgeStyle/.style = {thick}}
\tikzset{EdgeStyle/.append style = {->}}
\Edge(c0)(c1)
\Edge(c1)(c2)
\Edge(c2)(c3)
\Edge(c4)(c5)

\tikzset{EdgeStyle/.style = {bend right}}
\tikzset{EdgeStyle/.append style = {->}}
\Edge(c1)(c5)

\tikzset{EdgeStyle/.style = {bend left}}
\tikzset{EdgeStyle/.append style = {->}}
\Edge(c0)(c3)
\Edge[style={bend left=40}](c0)(c5)

%BOTTOM ROW

%Bottom Left
\Vertex[x=-3cm,y=-7.25cm,style={shade}]{d8}
\Vertex[x=-1cm,y=-7.25cm]{d2}
\Vertex[x=1cm,y=-7.25cm,style={shade}]{d9}
\Vertex[x=3cm,y=-7.25cm]{d3}
\Vertex[x=-3cm,y=-5.75cm]{d1}
\Vertex[x=-1cm,y=-5.75cm,style={shade}]{d11}
\Vertex[x=1cm,y=-5.75cm]{d0}
\Vertex[x=3cm,y=-5.75cm,style={shade}]{d6}
\Vertex[x=-3cm,y=-4.25cm]{d5}
\Vertex[x=-1cm,y=-4.25cm]{d4}
\Vertex[x=1cm,y=-4.25cm,style={shade}]{d10}
\Vertex[x=-1cm,y=-6.5cm,style={shade}]{d7}

\tikzset{EdgeStyle/.style = {thick}}
\Edge(d0)(d7)
\Edge(d1)(d8)
\Edge(d2)(d9)
\Edge(d4)(d11)
\Edge(d0)(d9)
\Edge(d1)(d11)

\tikzset{EdgeStyle/.style   = {line width = 2pt}}
\Edge(d0)(d6)
\Edge(d1)(d7)
\Edge(d2)(d8)
\Edge(d3)(d9)
\Edge(d4)(d10)
\Edge(d5)(d11)

\tikzset{EdgeStyle/.style   = {line width = 1pt}}
\draw (1.25cm,-6.5cm) node {{\Large $2$}};
\Edge[style={bend right=70}](d8)(d0)
\Edge[style={bend right=70}](d1)(d9)

%Bottom Middle
\Vertex[x=8.0cm,y=-7.0cm]{e0}
\Vertex[x=9.5cm,y=-7.0cm]{e1}
\Vertex[x=11cm,y=-7.0cm]{e2}
\Vertex[x=12.5cm,y=-7.0cm]{e3}
\Vertex[x=14cm,y=-7cm]{e4}
\Vertex[x=15.5cm,y=-7cm]{e5}
\Vertex[x=8.0cm,y=-4.5cm,style={shade}]{e6}
\Vertex[x=9.5cm,y=-4.5cm,style={shade}]{e7}
\Vertex[x=11cm,y=-4.5cm,style={shade}]{e8}
\Vertex[x=12.5cm,y=-4.5cm,style={shade}]{e9}
\Vertex[x=14cm,y=-4.5cm,style={shade}]{e10}
\Vertex[x=15.5cm,y=-4.5cm,style={shade}]{e11}

\Edge(e0)(e7)
\Edge(e1)(e8)
\Edge(e2)(e9)
\Edge(e4)(e11)
\Edge(e0)(e8)
\Edge(e0)(e9)
\Edge(e1)(e9)
\Edge(e1)(e11)

\draw (11.8cm,-4.6cm) node {{\Large $2$}};

\tikzset{EdgeStyle/.append style = {ultra thick}}
\Edge(e0)(e6)
\Edge(e1)(e7)
\Edge(e2)(e8)
\Edge(e3)(e9)
\Edge(e4)(e10)
\Edge(e5)(e11)

%Bottom Right
\Vertex[x=20.0cm,y=-5.75cm]{f0}
\Vertex[x=21.5cm,y=-5.75cm]{f1}
\Vertex[x=23cm,y=-5.75cm]{f2}
\Vertex[x=24.5cm,y=-5.75cm]{f3}
\Vertex[x=26cm,y=-5.75cm]{f4}
\Vertex[x=27.5cm,y=-5.75cm]{f5}

\tikzset{EdgeStyle/.style = {thick}}
\tikzset{EdgeStyle/.append style = {->}}
\Edge(f0)(f1)
\Edge(f1)(f2)
\Edge(f2)(f3)
\Edge(f4)(f5)

\tikzset{EdgeStyle/.style = {bend right}}
\tikzset{EdgeStyle/.append style = {->}}
\Edge[style={bend right=40}](f1)(f5)
\Edge(f1)(f3)

\tikzset{EdgeStyle/.style = {bend left}}
\tikzset{EdgeStyle/.append style = {->}}
\Edge(f0)(f2)
\Edge[style={bend left=40}](f0)(f3)

\draw (22.25cm,-4.5cm) node {{\Large $2$}};

%Arrows
\draw[->,thick](4.5,1.25) -- ( 6.5,1.25);
\draw[->,thick](16.5,1.25) -- (18.5,1.25);
\draw[->,thick](6.5,-5.75) -- ( 4.5,-5.75);
\draw[->,thick](18.5,-5.75) -- (16.5,-5.75);
\draw[->,thick](23.75,-.75) -- (23.75,-4.25);
\draw[->,thick,dashed](0,-1.25) -- (0,-3.25);
\end{tikzpicture}}
\end{center}
\caption{Constructing the parity closure of a bipartite graph with unique perfect matching}\label{finishedpicture}
\end{figure}
\end{Exam}

To emphasize, the constructed graph $G^+$ in the bottom-left is only
the \emph{potential} inverse of the original graph: There is no
guarantee that the eigenvalues of $G^+$ are the reciprocals of those
of $G$, but by Theorem \ref{gottabe}, \emph{if} the top-left graph is
invertible, then this is indeed the case.  Note then that the presence
of the double edge in the parity closure immediately implies that the
original graph is not simply invertible.  Regardless, the key issue
remaining is now to decide \emph{a priori} the invertibility of a
graph, a topic to which we devote the next section.  In particular, we
will see in Example \ref{promisedexample} that the graph in Figure \ref{finishedpicture}
is indeed invertible.

\medskip

Before doing so, let us close the current section by remarking on the
comment from the introduction that any reasonable notion of inversion
should have the property that the double-inverse of a graph should be
the original graph.  We prove that this phenomenon does indeed occur
with the parity closure.

\begin{Lemma}\label{G+isBUPM}
If $G$ is a bipartite graph with a unique perfect matching, then so is
$G^+$.  In particular, $G^{++}$ is defined.
\end{Lemma}
\begin{proof}
The block structure of $A^+$ provides the bipartition of $G^+$ and an
easy induction argument proves that the unique perfect matching of $G$
provides too the unique perfect matching of $G^+$.
\end{proof}

\begin{Theo}\label{G++}
If $G$ is an invertible bipartite graph with a unique perfect matching,
then $G^{++}=G$.
\end{Theo}
\begin{proof}
By Definition \ref{parclose}, the adjacency matrix of $G^+$ is
$|A^{-1}|$.  Since $G$ is invertible, there exists a signing matrix
$S$ such that the adjacency matrix $|A^{-1}|$ of $G^+$ can be
expressed $|A^{-1}|=SA^{-1}S$.  By Lemma \ref{G+isBUPM}, $G^{++}$
is defined, and we can consider its adjacency matrix
$||A^{-1}|^{-1}|=|(SA^{-1}S)^{-1}|=|SAS|$.  Since conjugation by $S$
only changes the signs of the entries of $A$, $|SAS|=|A|=A$, and we
have $G^{++}=G.$
\end{proof}

If we dub graphs satisfying the condition $G^{++}=G$ as
\emph{reflexive}, it is natural to wonder about the converse to the
theorem: Is every reflexive graph invertible?  We conclude this
section with a negative response to this question, via the reflexive
non-invertible counter-example below.  It would be interesting to
characterize invertible graphs among reflexive graphs.

\begin{figure}[!ht]\label{G++pic}
\begin{center}
\scalebox{.6}{
\begin{tikzpicture}
\GraphInit[vstyle=Normal]
\SetVertexMath
\SetVertexNoLabel

\Vertex[x=0.0cm,y=0.0cm]{a0}
\Vertex[x=2.0cm,y=0.0cm]{a1}
\Vertex[x=4.0cm,y=0.0cm]{a2}
\Vertex[x=6.0cm,y=0.0cm]{a3}
\Vertex[x=8.0cm,y=0.0cm]{a4}
\Vertex[x=10.0cm,y=0.0cm]{a5}

\AssignVertexLabel{a}{10}{$1$,$2$,$3$,$4$,$5$,$6$}
\tikzset{EdgeStyle/.append style = {->}}

\Edge(a0)(a1)
\Edge(a1)(a2)
\Edge(a2)(a3)
\Edge(a3)(a4)
\Edge[style={bend left=30}](a0)(a2)
\Edge[style={bend left=35}](a0)(a4)
\Edge[style={bend left=40}](a0)(a5)
\Edge[style={bend right=20}](a2)(a4)
\Edge[style={bend right=25}](a2)(a5)
\end{tikzpicture}}
\end{center}
\caption{A non-invertible reflexive graph.}
\end{figure}
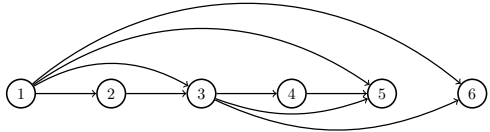

%*************************************************************

\section{Invertibility}

%*************************************************************

We maintain the notation from the previous section: $D$ is a directed
graph on $\{1,2,\ldots,n\}$ with adjacency matrix $B-I$.  Our goal is
to decide when $D$ is invertible, i.e., when $B^{-1}$ is signable to
its absolute value $|B^{-1}|,$ using the combinatorics of $D$.  The
following construction is of considerable interest.

\begin{Defn}
Given a digraph $D$, let $\Gamma=\Gamma_D$ denote its
\emph{maximal-path subgraph}, the spanning subgraph of $D$ which
includes an edge $e$ of $D$ if and only if there exist vertices $i$
and $j$ of $D$ such that $e$ lies on a path from $i$ to $j$ whose
length is maximal with respect to all such paths.
\end{Defn}

\begin{Remark}
Roughly, $\Gamma$ is the union, over all $i<j$, of all of the longest
paths from $i$ to $j$.  It will be occasionally useful to adopt the
equivalent but alternate point of view that $\Gamma$ is constructed by
deleting from $D$ every edge $e$ such that there exists a path in $D$
of length greater than one connecting the endpoints of $e$.  Thus we
will often speak of edges of $D$ ``surviving to $\Gamma$'' if no such
longer path exists.
\end{Remark}

The auxiliary graph $\Gamma$ is used significantly in \cite{TK09}
(where it is named the Hasse Diagram of $D$, consistent with the
poset-theoretic viewpoint therein) to prove the following necessary
condition for invertibility.  Our proof of the following result is
essentially different only in notation and terminology from the
original.

\begin{Theo}\label{GammaBip}
If $B$ is signable, then $\Gamma_D$ is bipartite.
\end{Theo}

\begin{proof}
It is trivial to reduce to the case that $\Gamma_D$ is connected, so
we assume this is the case.  Suppose there is a directed edge from $i$
to $j$ in $\Gamma_D$.  By definition of $\Gamma_D$, this forces all
directed paths from $i$ to $j$ in $D$ to have length one, and thus by
equation \eqref{Inv}, $B^{-1}_{i,j}<0$.  If $B^{-1}$ is signable, then
considering the entry $(SB^{-1}S)_{i,j}=s_iB^{-1}_{i,j}s_j\geq 0$, we
conclude that $s_i=-s_j,$ i.e., entries corresponding to vertices
connected by an edge of $\Gamma_D$ must have opposite signs.  By
connectedness, arbitrarily setting $s_1=+1$ (valid since $S$ signs
$B^{-1}$ if and only if $-S$ does) now induces the sign of $s_i$ for
each vertex $i$.  It is easy to see that grouping the vertices by this
sign provides a bipartition of the vertices of $\Gamma_D$.
\end{proof}
\begin{Defn}

We introduce a pairing $\<\cdot,\cdot\>\colon D\times D\to\{0,\pm1\}$
on the vertices of $D$ as follows: For $i,j\in D$, if $P$ is a
maximal-length path from $i$ to $j$ in $D$, let $\l(P)$ denote its
length, and set
\[
\<i,j\>=(-1)^{\l(P)}.
\]
If there are no paths from $i$ to $j$ in $D$, define $\<i,j\>=0$.
\end{Defn}
\begin{Lemma}\label{ijk}
If $\Gamma_D$ is bipartite and $i,j,k\in D$ are such that $i$ is
path-connected to $j$ and $j$ is path-connected to $k$, then
\[
\<i,j\>\<j,k\>=\<i,k\>.
\]
\end{Lemma}
\begin{proof}
By bipartiteness, a longest path from $i$ to $j$ concatenated with a
longest path from $j$ to $k$ must have the same parity as any (in
particular, the longest) path from $i$ to $k$.
\end{proof}
\noindent We now have the following necessary and sufficient condition for
invertibility.

\begin{Prop}\label{NecSuf}
If $\Gamma_D$ is bipartite, then $D$ is invertible if and only if for all $i,j\in
D$ we have
\begin{align}\label{necsuf}
\<i,j\>B_{i,j}^{-1}=\<i,j\>\hspace*{-.4cm}\sum_{P\in \P_{i,j}(D)}\hspace*{-.2cm}(-1)^{l(P)}\ge 0.
\end{align}
\end{Prop}
\begin{proof}
The equality is the definition of $B_{i,j}^{-1}$, so the content of
the proposition is the inequality.  We assume without loss of
generality that $D$ is connected.  For the first direction, assume
$\<i,j\>B_{i,j}^{-1}\geq 0$ for all $i,j$.  Choose a coloring of
$\Gamma_D$ and let $s_i=+1$ if vertex $i$ is colored black, and $-1$
otherwise.  Now $s_is_j$ is $+1$ if and only if $i$ and $j$ have the
same color, which since $\Gamma$ is bipartite and connected, occurs if
and only if each (undirected) path from $i$ to $j$ has even length.
Similarly, $s_is_j=-1$ if each path, in particular any
maximal-length path, has odd length.  In other words, we have
$s_is_j=\<i,j\>$.  Now, letting $S$ be the diagonal matrix whose
diagonal entries are the $s_i$, we have
\[
\<i,j\>B_{i,j}^{-1}\geq 0\quad\Longrightarrow \quad s_iB_{i,j}^{-1}s_j\geq 0,
\]
so $D$ is invertible by Theorem \ref{signable}.  Conversely, assume
that $D$ is invertible.  Then again by Theorem \ref{signable}, there
exists a signing matrix $S$ such that $SB^{-1}S=|B^{-1}|$, i.e.,
integers $s_i=\pm 1$ such that $s_iB_{i,j}^{-1}s_j> 0$ for all $i,j$
with $B_{i,j}^{-1}\neq 0$.  If $i$ and $j$ are adjacent in $\Gamma_D$,
then $\<i,j\>=-1$, and $B_{i,j}^{-1}<0$.  Thus $s_is_j=-1$, and
we have $s_is_j=\<i,j\>$ for all pairs and (\ref{necsuf}) is
satisfied.
\end{proof}

\begin{Remark}
Speaking loosely, the proposition tells us that a graph is invertible
if and only if the \emph{majority} of paths between any given pair of
vertices have the same parity as the \emph{longest} path between those
two vertices.  This is especially poignant, for example, when there
are few paths between each pair of vertices (Corollary
\ref{Bip2path}).
\end{Remark}

We note that the proposition is essentially equivalent to Corollary 5
in \cite{AK07}, which is proved using technical results from linear
algebra.  Our proof of this fact is self-contained and applies also to
multigraphs.  A motivating goal for the remainder of the section will
be to improve this result by honing in on the crucial substructures
which govern invertibility, in essence reducing the amount of
calculation needed to determine the invertibility of the graph.
Before doing so, let us extract from the proposition a few
corollaries.  The first of these is a generalization of Godsil's
Theorem 1 (\cite{Go85}) to multigraphs.

\begin{Cor}\label{GenGod}
If $D$ is bipartite, then $D$ is invertible and $D$ is a subgraph of
$D^+.$
\end{Cor}

\begin{proof}
For arbitrary vertices $i$ and $j$, the bipartiteness of $D$ ensures
that the lengths of all paths in $\P_{i,j}(D)$ have the same parity.
Thus, every term in the sum $\sum_{P\in\P_{i,j}(D)}(-1)^{l(P)}$ has
the same sign, and the same sign as $\<i,j\>$.  Therefore, by Theorem
\ref{NecSuf}, $D$ is invertible.  Furthermore, $|B_{i,j}^{-1}|$ equals
the number of paths between $i$ and $j$, and thus $B_{i,j}\le
|B_{i,j}^{-1}|$, i.e., that $D$ is a subgraph of $D^+$.
\end{proof}

\begin{Cor}
If $D$ is a simple, invertible graph and $D^+$ is bipartite, then $D$
is simply invertible.
\end{Cor}

\begin{proof}
Combining Theorem \ref{G++} and Corollary \ref{GenGod}, $D^+$ is a
subgraph of the simple graph $D^{++}=D$.
\end{proof}

\begin{Lemma}\label{intersection}
For any $D$, $\Gamma_D$ is also a subgraph of $D^+$ (and hence of
$D\cap D^+$).
\end{Lemma}
\begin{proof}
Suppose $e$ is an edge from $i$ to $j$ in $\Gamma_D$.  We prove that
$B_{i,j}^{-1}\neq 0$, showing that there exists an edge from $i$ to
$j$ in $D^+$.  Indeed, by definition of $B_{i,j}^{-1}$, the existence
of an odd-length path from $i$ to $j$ in $D$ (the path consisting of
$e$ alone) implies that if $B_{i,j}^{-1}$ were equal to zero, there
would also be at least one even-length path from $i$ to $j$ in $D$.
This contradicts that $e$ lied on a maximal-length path in $D$.
\end{proof}

\begin{Theo}\label{BipEquiv}
For any $D$, the following are equivalent:
\begin{itemize}
\item[(i)] $D$ and $D^+$ are both bipartite.
\item[(ii)] There is no path of length greater than one in $D$.
\item[(iii)] $D=D^+$ and $D$ is bipartite.
\end{itemize}
\end{Theo}

\begin{proof}
Clearly (i) follows directly from (iii).  It remains to prove
$(\text{i}) \Rightarrow (\text{ii})\Rightarrow (\text{iii})$.  First,
$(\text{i})\Rightarrow (\text{ii}):$ Suppose $D$ and $D^+$ are bipartite and that
there exists a path of length greater than one in $D$.  Such a path
ensures that there exist vertices $i$ and $j$ connected by an even
length path in $D$.  We conclude from this two consequences:
\begin{itemize}
\item The bipartiteness of $D$ implies there is no odd path between
  $i$ and $j$ in $D$.  Thus, by (\ref{necsuf}),
  $B^{-1}_{i,j}=\sum_{P\in\P_{i,j}(D)}(-1)^{l(P)}>0$ and so vertices
  $i$ and $j$ are adjacent in $D^+$.
\item By definition of $\Gamma_D$ and the bipartiteness of $D$, since
  there is an even length path between vertices $i$ and $j$ in $D$,
  then there is also an even length path between $i$ and $j$ in
  $\Gamma_D$.  Since $\Gamma_D$ is a subgraph of $D^+$ (Lemma
  \ref{intersection}), there is an even length path between adjacent
  vertices $i$ and $j$ in $D^+$.
\end{itemize}
These two observations contradict that $D^+$ was assumed bipartite.
Finally, we prove $(\text{ii})\Rightarrow (\text{iii}):$ Suppose there
does not exist a path of length greater than one in $D$.  From the
definition of parity closure, for any two vertices $i$ and $j$, the
number of edges from vertex $i$ to vertex $j$ in $D^+$
is $$|\#(\mbox{even paths from $i$ to $j$ in $D$})-\#(\mbox{odd paths
  from $i$ to $j$ in $D$})|.$$ Since each such path is length 1, this
simply returns the number of edges from $i$ to $j$ in $D$.  Thus,
$D=D^+$.  Furthermore, the assumption that the length of all paths in
$D$ is less than or equal to 1 implies $D$ is bipartite.
\end{proof}

Our next corollary is the extension of the main theorem of Simion-Cao
\cite{SC89} to multigraphs.  Here, we define a graph to be
\emph{self-dual} if $G=G^+$ and the \emph{corona} of a graph $H$ to be
the graph obtained by adding to $H$ a neighbor of degree 1 to each
vertex.

\begin{Cor}\label{SimionCao}
A graph $G$ is the corona of a bipartite multigraph $H$ if and only if
$G$ is a self-dual bipartite multigraph with a unique perfect matching
and bipartite $D$.
\end{Cor}

\begin{proof}
Suppose $G$ is the corona of a bipartite multigraph $H$.  Clearly $G$
inherits the bipartiteness of $H$ and contains the unique perfect
matching $M$ consisting of the $n$ pendant edges.  Furthermore, since
each edge of $M$ contains a pendant vertex, each non-isolated vertex
of $D$ is either a source vertex or a sink vertex (i.e., can have
vertices adjacent to it or is adjacent to other vertices, but not
both).  Thus $D$ contains no paths of length greater than one.  By
Theorem \ref{BipEquiv}, $D=D^+$ and $D$ is bipartite.

\medskip

The converse follows similarly:  Suppose $G$ is a bipartite multigraph
with unique perfect matching $M$, $D$ bipartite, and $G=G^+.$ By
Theorem \ref{BipEquiv}, $D$ does not contain a path of length
greater than one.  Thus each vertex of $D$ is at most either a
source or a sink (in particular, not both), and hence each edge in $M$ contains
a pendant vertex.  Thus $G$ is the corona of $D$.
\end{proof}

%*************************************************************

\subsection{Algorithm for Determining Invertibility}\label{algorithm-subsec}

%*************************************************************

\medskip

We begin with some combinatorial preliminaries on path-sets.  For a
graph $D$ on vertices labelled $\{1,2,\ldots,n\}$ and a vertex $1\leq
k\leq n$, let $D-\{k\}$ be the subgraph of $D$ resulting from the
deletion of vertex $k$, retaining the original labels on all other
vertices.  Recall that for a digraph $D$, we denote by $B=B_D$ the
upper uni-triangular matrix such that $B-I$ is the adjacency matrix of
$D$.  The following lemma relates the path-counting in $D$ to that of
$D-\{k\}$.

\begin{Lemma}\label{algLem}
Suppose $i<k<j$ and let ${_k}{B_{i,j}^{-1}}$ denote the $(i,j)^{th}$
entry of $B_{D-\{k\}}^{-1}.$ Then
$$B_{i,j}^{-1}=B_{i,k}^{-1}B_{k,j}^{-1}+{_k}{B_{i,j}^{-1}}.$$
\end{Lemma}

\begin{proof}
We abbreviate $\P_{i,j}=\P_{i,j}(D)$.  If ${_k}\P_{i,j}$ denotes the
set of paths from $i$ to $j$ \emph{not} passing through $k$, then we
have the decomposition $\P_{i,j}=(\P_{i,k}\times \P_{k,j})\cup
{_k}\P_{i,j}$, where we interpret an element of the product of the two
path-sets as the concatenation of the two paths.  Now by Equation
\eqref{Bij}, we have
\begin{align*}
B_{i,j}^{-1}=\sum_{P\in\P_{i,j}}(-1)^{l(P)}&=\sum_{P\in\P_{i,k}}\sum_{P'\in\P_{k,j}}(-1)^{l(P)+l(P')}+\sum_{P\in{_k}\P_{i,j}}(-1)^{l(P)}\\
&=\sum_{P\in\P_{i,k}}(-1)^{l(P)}\sum_{P'\in\P_{k,j}}(-1)^{l(P')}+\sum_{P\in \P_{i,j}(D-\{k\})}(-1)^{l(P)}\\
&=B_{i,k}^{-1}B_{k,j}^{-1}+{_k}B_{i,j}^{-1}.
\end{align*}

\end{proof}

We remark that the product structure induced from concatenation of
paths is the key simplifying tool: a pair $(i,j)$ has the property
referenced at the start of the section -- that there exists an
intermediate vertex $k$ such that every path from $i$ to $j$ passes
through $k$ -- if and only if the set $\P_{i,j}$ can be written
(``factors'') as the product of non-empty path-sets $\P_{i,k}$ and
$\P_{k,j}$.  This loose analogy between factoring integers and
path-sets motivates the following terminology which will be of some use.

\begin{Defn}
Given vertices $i$ and $j$ in $D$, we call the ordered pair $(i,j)$:
\begin{itemize}
\item a \emph{zero pair} if there is no path in $D$ from $i$ to $j$;
\item a \emph{unit pair} if all paths from $i$ to $j$ in $D$ have length 1;
\item a \emph{composite pair} if there exists a vertex $k\in D$ with
  $i<k<j$ such that every path $P$ from $i$ to $j$ in $D$ passes through $k$;
\item a \emph{prime pair} if it is neither a zero pair, a unit pair,
  nor a composite pair.
\end{itemize}
\end{Defn}

Finally, for graphs with $\Gamma_D$ bipartite, we say the pair $(i,j)$
is \emph{signable} if $\<i,j\>B_{i,j}^{-1}\ge 0$ and \emph{simply
  signable} if $\<i,j\>B_{i,j}^{-1}\in\{0,1\}$.  Note that under this
definition of the signability of a pair, Proposition \ref{NecSuf} says
that a graph is invertible (resp. simply invertible) if and only if
\emph{all} pairs $(i,j)$ are signable (resp. simply signable).

\smallskip

We turn to streamlining the results of Proposition \ref{NecSuf}.  In
particular, we look to significantly reduce the number of pairs
$(i,j)$ whose signability we need to evaluate in order to determine
the invertibility of the graph.  The rough strategy is to note that by
Lemma \ref{algLem}, if all paths from $i$ to $j$ pass through an
intermediate vertex $k$, then the signability of $(i,j)$ is forced by
the signability of the pairs $(i,k)$ and $(k,j)$.

\medskip

\begin{Theo}\label{Alg}
A digraph $D$ is invertible (resp. simply invertible) if and only if
$\Gamma_D$ is bipartite and all prime pairs of $D$ are signable.
Similarly, a digraph $D$ is simply invertible if and only if
$\Gamma_D$ is simple and bipartite and all prime pairs of $D$ are
simply signable.
\end{Theo}

\begin{proof}
If $D$ is invertible, then by Theorem \ref{NecSuf} and Lemma
\ref{GammaBip}, $\Gamma_D$ is bipartite and
\[
\<i,j\>\hspace*{-.4cm}\sum_{P\in \P_{i,j}(D)}\hspace*{-.2cm}(-1)^{l(P)}=\<i,j\>B_{i,j}^{-1}\ge 0
\]
for all prime pairs $(i,j).$ For the converse, suppose that $\Gamma_D$
is bipartite and that all prime pairs are signable, and let $(i,j)$ be
an arbitrary non-prime pair.  If $(i,j)$ is a zero pair, it is
trivially signable, and if $(i,j)$ is a unit pair, then $\<i,j\>=-1$
and $B_{i,j}^{-1}$ is negative the number of edges from $i$ to $j$,
giving $\<i,j\>B_{i,j}^{-1}>0$.  Finally, assume $(i,j)$ is composite,
so there exists a $k$ with $i<k<j$ and $\P_{i,j}=\P_{i,k}\P_{k,j}$.
Since no paths from $i$ to $j$ omit $k$, we have by Lemma \ref{algLem}
that $B_{i,j}^{-1}=B_{i,k}^{-1}B_{k,j}^{-1}+0$.  Recalling that prime
pairs are signable, by induction we can assume that
$\<i,k\>B_{i,k}^{-1}\ge 0$ and $\<k,j\>B_{k,j}^{-1}\ge 0$, and so
\begin{equation}\label{compositepair}
\<i,j\>B_{i,j}^{-1}=\<i,k\>B_{i,k}^{-1}\,\<k,j\>B_{k,j}^{-1}\ge
0,
\end{equation}
proving that in fact all pairs of $D$ are signable.  The claim for
simple invertibility proceeds similarly: $D$ is clearly simply
invertible if and only if \emph{all} pairs $(i,j)$ are simply
signable.  It remains to show that it is sufficient to check the prime
pairs.  Since $\Gamma_D$ is simple, all unit pairs are simply
signable, and all zero pairs are trivially simply signable.  If we
assume all prime pairs are simply signable, then equation
\eqref{compositepair} shows that all composite pairs (and so all
pairs) are simply signable, as desired.
\end{proof}

\noindent This strengthening of Proposition \ref{NecSuf} immediately decides
the invertibility of some large classes of graphs.

\begin{Cor}\label{Bip2path}
If $\Gamma_D$ is bipartite and if there are two or fewer paths between
each prime pair $(i,j)$ of $D$, then $D$ is invertible.
\end{Cor}
\begin{proof}
Note by definition of $(i,j)$ being a prime pair, there cannot be zero
or one paths from $i$ to $j$.  If there are two paths $P$ and $P'$
with, say, $\l(P)\geq \l(P')$, then
\[
\<i,j\>B_{i,j}^{-1}=(-1)^{\l(P)}((-1)^{\l(P)}+(-1)^{\l(P')})=1+(-1)^{\l(P)+\l(P')}\geq 0.
\]
Hence all prime pairs are signable.
\end{proof}

\noindent We include next a partial converse of Godsil's Theorem:
\begin{Cor}
A digraph $D$ without prime pairs is invertible if and only if
it is bipartite.
\end{Cor}
\begin{proof}
By Theorem \ref{Alg}, $D$ is invertible if and only if $\Gamma_D$ is
bipartite, but if $D$ has no prime pairs, then $D=\Gamma_D$.
\end{proof}

Before turning to some more involved examples, we observe that
determining invertibility via even the reduced process of checking
only prime pairs involves some redundant calculations.  Specifically,
since there can be substantial overlap in the computations needed for
verifying the signability of a prime pair, we can avoid (or at least
attempt to minimize) redundancy by attending first to prime pairs
$(i,j)$ with smaller values of $|j-i|$.  Since we will wish to address
the secondary question of when an invertible graph is additionally
\emph{simply} invertible, we note that in light of the second claim of
Theorem \ref{Alg} it is sufficient to check that $|B_{i,j}^{-1}|\le 1$
for all prime pairs $(i,j)$.  Recall that ${_k}B_{i,j}^{-1}$ denotes
the $(i,j)^{th}$ entry of $B_{D-\{k\}}^{-1}.$

\begin{Exam}\label{introleft}
This example determines the invertibility of the graph $G$ in Figure \ref{TKfigs}
(left) from the introduction.  Presented is the associated graph $D$, with
its subgraph $\Gamma_D$ consisting of the bold edges.

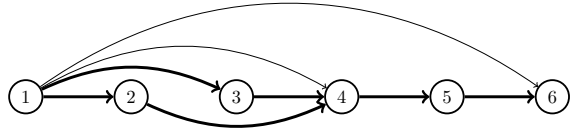
\begin{figure}[!ht]
\begin{center}\scalebox{.7}{
\begin{tikzpicture}
\GraphInit[vstyle=Normal]
\SetVertexMath
\SetVertexNoLabel

\Vertex[x=0.0cm,y=0.0cm]{a0}
\Vertex[x=2.0cm,y=0.0cm]{a1}
\Vertex[x=4.0cm,y=0.0cm]{a2}
\Vertex[x=6.0cm,y=0.0cm]{a3}
\Vertex[x=8.0cm,y=0.0cm]{a4}
\Vertex[x=10.0cm,y=0.0cm]{a5}

\AssignVertexLabel{a}{10}{$1$,$2$,$3$,$4$,$5$,$6$}
\tikzset{EdgeStyle/.append style = {->,ultra thick}}
\Edge(a0)(a1)
\Edge(a2)(a3)
\Edge(a3)(a4)
\Edge(a4)(a5)
\tikzset{EdgeStyle/.style = {->,bend left,ultra thick}}
\Edge[style={bend left=25}](a0)(a2)
\tikzset{EdgeStyle/.style = {->,bend right,ultra thick}}
\Edge[style={bend left=25}](a1)(a3)
\tikzset{EdgeStyle/.style = {->,bend left, thin}}
\Edge[style={bend left=30}](a0)(a3)
\Edge[style={bend left=35}](a0)(a5)

\end{tikzpicture}}
\end{center}
\end{figure}
\noindent Since $\Gamma_D$ contains no odd cycles, it is bipartite.
The only prime pairs of $D$ are $(1,4)$ and
$(1,6).$ We have
\[
\<1,4\>B_{1,4}^{-1}=(-1)^2\left(\#(\mbox{even paths from $1$ to
  $4$})-\#(\mbox{odd paths from $1$ to $4$})\right)=2-1=1\geq 0,\] so
$(1,4)$ is signable.  By Lemma
\ref{algLem}, $$B_{1,6}^{-1}=B_{1,4}^{-1}B_{4,6}^{-1}+{_4}B_{1,6}^{-1}=(1)(1)+(-1)=0,$$
and so $(1,6)$ is also signable.  Further, since $|B_{i,j}^{-1}|\leq
1$ for both prime pairs, $D$ (and hence $G$) is simply invertible.
\end{Exam}

\begin{Exam}\label{introright}
Similarly, below is the graph $D$ for the graph $G$ in Figure \ref{TKfigs} (right) of
the introduction.  Again, the maximal-path subgraph $\Gamma$ is drawn in bold.
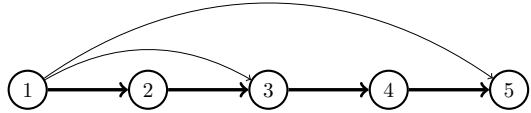
\begin{center}
\scalebox{.8}{
\begin{tikzpicture}
\GraphInit[vstyle=Normal]
\SetVertexMath
\SetVertexNoLabel

\Vertex[x=0.0cm,y=0.0cm]{a0}
\Vertex[x=2.0cm,y=0.0cm]{a1}
\Vertex[x=4.0cm,y=0.0cm]{a2}
\Vertex[x=6.0cm,y=0.0cm]{a3}
\Vertex[x=8.0cm,y=0.0cm]{a4}

\AssignVertexLabel{a}{10}{$1$,$2$,$3$,$4$,$5$}
\tikzset{EdgeStyle/.append style = {->,ultra thick}}
\Edge(a0)(a1)
\Edge(a1)(a2)
\Edge(a2)(a3)
\Edge(a3)(a4)
\tikzset{EdgeStyle/.style = {thin,->,bend left=30}}
\Edge(a0)(a2)
\tikzset{EdgeStyle/.style = {thin,->,bend left=35}}
\Edge(a0)(a4)

\end{tikzpicture}}
\end{center}
First we verify $\Gamma_D$ is bipartite and note the only prime pairs
are $(1,3)$ and $(1,5).$  We have $B_{1,3}^{-1}=1-1=0$, but
\[
\<1,5\>B_{1,5}^{-1}=(-1)^4(B_{1,3}^{-1}B_{3,5}^{-1}+{_3}B_{1,5}^{-1})=(0)(1)+(-1)=-1,
\]
so $(1,5)$ is not signable, and $G$ is not invertible.
\end{Exam}

\begin{Exam}\label{promisedexample}
Finally, we reconsider the example of Figure \ref{finishedpicture}, whose invertibility
was left unanswered.
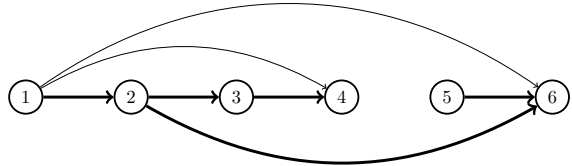
\begin{center}
\scalebox{.7}{
\begin{tikzpicture}
\GraphInit[vstyle=Normal]
\SetVertexMath
\SetVertexNoLabel

\Vertex[x=0cm,y=0cm]{c0}
\Vertex[x=2cm,y=0cm]{c1}
\Vertex[x=4cm,y=0cm]{c2}
\Vertex[x=6cm,y=0cm]{c3}
\Vertex[x=8cm,y=0cm]{c4}
\Vertex[x=10cm,y=0cm]{c5}

\AssignVertexLabel{c}{10}{$1$,$2$,$3$,$4$,$5$,$6$}
\tikzset{EdgeStyle/.append style = {->,ultra thick}}
\Edge(c0)(c1)
\Edge(c1)(c2)
\Edge(c2)(c3)
\Edge(c4)(c5)

\tikzset{EdgeStyle/.style = {ultra thick,->,bend right=30}}
\Edge(c1)(c5)
\tikzset{EdgeStyle/.style = {thin,->,bend left=30}}
\Edge(c0)(c3)
\tikzset{EdgeStyle/.style = {thin,->,bend left=35}}
\Edge(c0)(c5)
\end{tikzpicture}}
\end{center}

\noindent The only prime pairs are $(1,4)$ and $(1,6)$.  By counting
paths, we have
\[
\<1,6\>B_{1,6}^{-1}=(-1)^2(1+(-1))=0
\]
and
\[
\<1,4\>B_{1,4}^{-1}=(-1)^3(0-2)=2\geq 0,
\]
so both prime pairs are signable.  Hence $G$ is invertible, but since
$|B_{1,4}^{-1}|>1$, it is not simply invertible.
\end{Exam}

%%%%%%%%%%%%%%%%%%%%%%%%%%%%%%%%%%%%%%%%%%%%%%%%%%%%%%%%%%%%%%%%%%%%

\subsection{Unicyclic graphs}\label{UInvert}

Several authors have examined the invertibility of unicyclic bipartite
graphs with a unique perfect matching.  In \cite{AK07}, Akbari and
Kirkland present necessary and sufficient criteria for the
invertibility of such graphs.  Below, we establish some preliminaries
and derive their criteria from a prime pairs perspective.  Briefly,
the point is that a unicyclic graph can have at most one prime pair,
so the methods from the last section are particularly apt.  In Section
\ref{unicycleconstruct}, we present a new algorithm for constructing
\emph{all} invertible unicyclic bipartite graphs.

\medskip

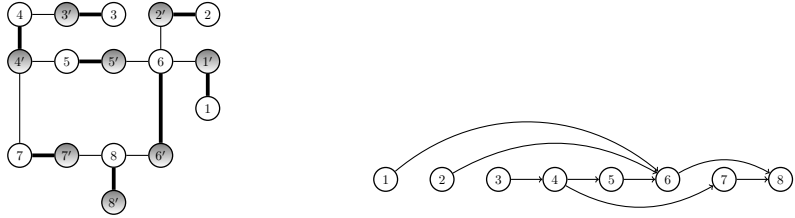
\begin{figure}[!ht]
\begin{center}
\begin{tabular}{ccc}
\scalebox{.5}{
\begin{tikzpicture}
\GraphInit[vstyle=Normal]
\SetVertexMath
\SetVertexNoLabel

\Vertex[x=0.0cm,y=0.0cm]{a0}
\Vertex[x=0.0cm,y=-1.25cm,style={shade}]{a1}
\Vertex[x=0.0cm,y=-3.75cm]{a2}
\Vertex[x=1.25cm,y=0.0cm,style={shade}]{a3}
\Vertex[x=1.25cm,y=-1.25cm]{a4}
\Vertex[x=1.25cm,y=-3.75cm,style={shade}]{a5}
\Vertex[x=2.5cm,y=0.0cm]{a6}
\Vertex[x=2.5cm,y=-1.25cm,style={shade}]{a7}
\Vertex[x=2.5cm,y=-3.75cm]{a8}
\Vertex[x=2.5cm,y=-5.0cm,style={shade}]{a9}
\Vertex[x=3.75cm,y=0.0cm,style={shade}]{a10}
\Vertex[x=3.75cm,y=-1.25cm]{a11}
\Vertex[x=3.75cm,y=-3.75cm,style={shade}]{a12}
\Vertex[x=5.0cm,y=0.0cm]{a13}
\Vertex[x=5.0cm,y=-1.25cm,style={shade}]{a14}
\Vertex[x=5.0cm,y=-2.5cm]{a15}

\AssignVertexLabel{a}{10}{$4$,$4'$,$7$,$3'$,$5$,$7'$,$3$,$5'$,$8$,$8'$,$2'$,$6$,$6'$,$2$,$1'$,$1$}

\Edge(a0)(a3)
\Edge(a1)(a4)
\Edge(a1)(a2)
\Edge(a7)(a11)
\Edge(a10)(a11)
\Edge(a11)(a14)
\Edge(a5)(a8)
\Edge(a8)(a12)
\renewcommand*{\EdgeLineWidth}{1mm}
\tikzset{EdgeStyle/.style   = {line width = \EdgeLineWidth}}
\Edge(a3)(a6)
\Edge(a10)(a13)
\Edge(a0)(a1)
\Edge(a4)(a7)
\Edge(a2)(a5)
\Edge(a8)(a9)
\Edge(a11)(a12)
\Edge(a14)(a15)

\end{tikzpicture}}
&
\hspace*{.5in}
&

\scalebox{.5}{
\begin{tikzpicture}
\GraphInit[vstyle=Normal]
\SetVertexMath
\SetVertexNoLabel

\Vertex[x=7.0cm,y=1.5cm]{a0}
\Vertex[x=8.5cm,y=1.5cm]{a1}
\Vertex[x=10.0cm,y=1.5cm]{a2}
\Vertex[x=11.5cm,y=1.5cm]{a3}
\Vertex[x=13.0cm,y=1.5cm]{a4}
\Vertex[x=14.5cm,y=1.5cm]{a5}
\Vertex[x=16.0cm,y=1.5cm]{a6}
\Vertex[x=17.5cm,y=1.5cm]{a7}

\AssignVertexLabel{a}{10}{$1$,$2$,$3$,$4$,$5$,$6$,$7$,$8$}
\tikzset{EdgeStyle/.append style = {->}}
\Edge(a2)(a3)
\Edge(a3)(a4)
\Edge(a4)(a5)
\Edge(a6)(a7)
\tikzset{EdgeStyle/.style = {->,bend left}}
\Edge(a1)(a5)
\Edge(a5)(a7)
\Edge[style={bend left=40}](a0)(a5)
\tikzset{EdgeStyle/.style = {->,bend right}}
\Edge(a3)(a6)

\end{tikzpicture}}
\end{tabular}
\end{center}
\caption{A bipartite unicyclic graph and its associated digraph.}\label{unicycle}
\end{figure}

Let $U$ be a unicyclic bipartite graph with unique perfect matching.
The bipartiteness of $U$ ensures that the cycle has even length, while
the unique perfect matching forces the number of edges which are both
in the perfect matching and incident to (and not in) the cycle to be
even.  Let $2m$ be the length of the cycle and $2k$ be the number of
matched edges incident to the cycle.  The topological sort of the
vertices in Section \ref{prelims} can be chosen to also consecutively
order the matched edges in or incident to the cycle.  The associated
digraph $D$ thus contains a single undirected cycle of length $m+k$,
say on consecutive vertices $\{v+1,\ldots, v+(m+k)\}$, corresponding
to the matched edges in or incident to the unicycle.  Because a prime
pair in any digraph produces an undirected cycle, the only possible
prime pair in $D$ is $(v+1,v+m+k)$.  Figure \ref{unicycle} shows a bipartite
unicyclic graph with unique perfect matching and its associated
digraph, with $m=4$ and $k=1$.  Note that vertices $\{4,5,6,7,8\}$
form the undirected cycle of $D$.

\noindent In the following, we abuse terminology slightly and say that
$D$ is unicyclic if it is unicyclic as an undirected graph.

\begin{Lemma}\label{ppk}
The pair $(v+1,v+m+k)$ is prime if and only if $k=1$.
\end{Lemma}

\begin{proof}
Let $i=v+1$ and $j=v+m+k$.  Consider the directed subgraph $D'$ of $D$
induced by the vertices $\{i,\ldots,j\}$ corresponding to the matched
edges in or incident to the unicycle in $U$.  For all
$w\in\{i,\ldots,j\}$, $w$ is a source or sink vertex in $D'$ if and
only if $w$ corresponds to an incident matched edge.  Thus, the number
of source/sink vertices in $D'$ is $2k$.  The pair $(i,j)$ forms a
prime pair in $D$ if and only if $i$ is the only source and $j$ is the
only sink in $D'$, i.e., if and only if $k=1$.
\end{proof}

\begin{Theo}\label{unicycleinvert}
A unicyclic digraph $D$ is invertible if and only if $m+k$ is even or $k=1$
and the two matched edges incident to the unicycle are incident to
adjacent vertices.
\end{Theo}

\begin{proof}

Again let $i=v+1$ and $j=v+m+k$, so that $\{i,\ldots,j\}$ is the set
of vertices forming an undirected cycle in $D$ (listed in increasing
order).  Now $\Gamma_{D}=D$ unless $(i,j)$ forms a prime pair and $i$
is adjacent to $j$.  Thus, by Lemma \ref{ppk}, if $k\not=1$ or $i$ is
not adjacent to $j$, and $m+k$ is odd, then $\Gamma_{D}=D$ and
contains an odd cycle of length $m+k$.  By Theorem \ref{GammaBip}, $D$
is not invertible.

\smallskip

If $m+k$ is even, or $k=1$ and vertices $i$ and $j$ are adjacent, then
$\Gamma_{D}$ either contains only an even cycle or $\Gamma_{D}$ is the
acyclic graph formed by deleting the edge between $i$ and $j$ in $D$.
In either case, $\Gamma_{D}$ is bipartite.  The unicyclicity of $D$
ensures there are two or fewer paths between the only possible prime
pair $(i,j)$.  Thus, $D$ is invertible by Corollary \ref{Bip2path}.
\end{proof}

\begin{Cor}
An invertible unicyclic graph $D$ is \emph{simply} invertible if and
only if $k>1$ or $m+k$ is odd.
\end{Cor}

\begin{proof}
If $k=1$ then there are exactly two paths between the only prime pair
vertices $i$ and $j$.  If those paths are of the same parity (i.e., if
$m+k$ is even) then $|B_{i,j}^{-1}|=2$ and $D$ is not simply
invertible.  If those paths are of opposite parity (i.e., $m+k$ is
odd) then $|B_{i,j}^{-1}|=0$ and $D$ is simply invertible.  If $k>1$
then $D$ contains no prime pairs and $|B_{x,y}^{-1}|\le 1$ for all
vertices $x, y$ in $D$.
\end{proof}

%%%%%%%%%%%%%%%%%%%%%%%%%%%%%%%%%%%%%%%%%%%%%%%%%%%%%%%%%%%%%%%%%%%%%%%%%%%%%%%%%%

\subsection{Optimal Bipartite Subgraph}

We finish this section with an attempt to unite the two striking links
between invertibility and bipartiteness found in \cite{Go85} and
\cite{TK09}, and their generalization to multigraphs (Corollary
\ref{GenGod} and Theorem \ref{GammaBip}).  Namely, we have that for a
directed graph $D$:

\begin{itemize}
\item If $D$ is bipartite, $D$ is invertible.
\item If $D$ is invertible, $\Gamma_D$ is bipartite.
\end{itemize}
\noindent
One is led to wonder the extent to which these two results are optimal, i.e., ask:
\begin{itemize}
\item Are there proper subgraphs of $D$ whose bipartiteness implies
  the invertibility of $D$?
\item Are there proper supergraphs of $\Gamma_D$ whose bipartiteness
  is implied by the invertibility of $D$?
\end{itemize}
The answers to these can be trivially negative for a \emph{given}
graph (e.g., take any invertible graph such that $D=\Gamma_D$), so we
mean to ask these in the broader context of canonical constructions
that have the potential of being proper subgraphs and supergraphs,
respectively.  In particular, it is natural to ask if there is a
canonically-defined subgraph $H_D$ of $D$ satisfying $\Gamma_D\leq
H_D\leq D$ whose bipartiteness is equivalent to the invertibility of
$D$.  In this context, we answer both of the bulleted question in the
affirmative, but conclude that such a ``bipartite litmus test'' does
not exist.  Let us first improve the ``upper bound'' $H_D\leq D$.

\begin{Theo}\label{minimalprimepairs}
Suppose $(i_1,j_1), (i_2,j_2),\ldots,(i_n,j_n)$ are prime pairs of $D$
such that for all $k$:
\begin{itemize}
\item $j_k\leq i_{k+1}$ for $1\leq k\leq n-1$.
\item There is an edge $e_k$ from $i_k$ to $j_k$ in $D$.
\end{itemize}
Let $D'$ be the graph obtained from $D$ by deleting $e_k$, for $1\leq
k\leq n$.  Then if $D'$ is bipartite, $D$ is invertible.
\end{Theo}
\begin{proof}
By assumption, the removed edges $e_k$ did not survive to $\Gamma_D$,
so $\Gamma_D$ is a subset of $D'$.  In particular, since $D'$ is
assumed bipartite, $\Gamma_D$ is as well, and so we need only to check
that each prime pair $(a,b)$ of $D$ is signable.  We partition the set
$\P=\P_{a,b}(D)$ of paths from $a$ to $b$ in $D$ according to which of
the pairs $(i_k,j_k)$ a path goes through, i.e., writing
\[
\P=\bigcup_{S\subset\{1,2,\ldots,n\}}\P_S,
\]
where $\P_S$ is the subset of $\P$ consisting of paths that pass
through both $i_k$ and $j_k$ for exactly those $k\in S$ (of course,
$S$ can be empty).  This induces an analogous decomposition of
$B_{a,b}^{-1}$:
\[
\<a,b\>B_{a,b}^{-1}=\<a,b\>\sum_{P\in \P}(-1)^{\l(P)}=\sum_{S\subset\{1,2,\ldots,n\}}\<a,b\>\sum_{P\in \P_S}(-1)^{\l(P)}.
\]
It suffices to show that each of these summands is non-negative.
Write $S=\{m_1,\ldots,m_k\}$, so that
\[
\P_S=\P_{a,i_{m_1}}\times\P_{i_{m_1},j_{m_1}}\times\P_{j_{m_1},i_{m_2}}\times\P_{i_{m_2},j_{m_2}}\times\cdots \times\P_{j_{m_{k-1}},i_{m_k}}\times\P_{i_{m_k},j_{m_k}}\times\P_{j_{m_k},b}.
\]
Setting $j_{m_0}=a$ and $i_{m_{k+1}}=b$ for notational
convenience (if not aesthetics), we find
\[
\<a,b\>\sum_{P\in \P_S}(-1)^{\l(P)}=\prod_{\alpha=0}^k\<j_{m_\alpha},i_{m_{\alpha+1}}\>B_{j_{m_\alpha},i_{\alpha+1}}^{-1}\prod_{\alpha=1}^k\<i_{m_\alpha},j_{m_\alpha}\>B_{i_{m_\alpha},j_{m_\alpha}}^{-1}.
\]
Here we have used Lemma \ref{ijk}.  Finally, we check that each of the
factors above is positive:
\begin{enumerate}[(1)]
\item Since the path-sets $\P_{j_{m_\alpha},i_{m_\alpha+1}}$ remain
  unchanged by the deletion of the $e_k$'s, the corresponding factors
  can be written as follows:
\[\<j_{m_\alpha},i_{m_{\alpha+1}}\>B_{j_{m_\alpha},i_{\alpha+1}}^{-1}(D)=\<j_{m_\alpha},i_{m_{\alpha+1}}\>B_{j_{m_\alpha},i_{\alpha+1}}^{-1}(D').\]
Now each of these are positive since $D'$ is bipartite, making every
pair $(j_{m_\alpha},i_{m_\alpha+1})$ signable.
\item Since the path sets $\P_{i_{m_\alpha},j_{m_\alpha}}(D)$ and
  $\P_{i_{m_\alpha},j_{m_\alpha}}(D')$ differ by the single edge
  $e_{m_\alpha}$, the remaining factors satisfy
\[
\<i_{m_\alpha},j_{m_\alpha}\>B_{i_{m_\alpha},j_{m_\alpha}}^{-1}=\<i_{m_\alpha},j_{m_\alpha}\>B_{i_{m_\alpha},j_{m_\alpha}}^{-1}(D')\pm 1\geq 0,
\]
where we have used that by virtue of signability and having at least
one other path between them, $\<i_{m_\alpha},j_{m_\alpha}\>B_{i_{m_\alpha},j_{m_\alpha}}^{-1}(D')$
is positive and at least one.
\end{enumerate}
\end{proof}
\begin{Remark}\label{proper}
We note that if $D\neq \Gamma_D$, then one can find prime pairs as in
the theorem, thereby constructing a \emph{proper} subgraph $D'$ of $D$
whose bipartiteness implies the invertibility of $D$.
\end{Remark}

\begin{Exam}
The figure below gives a (non-bipartite) digraph the demonstration of
whose invertibility requires a moderate calculation using the
techniques of Section \ref{algorithm-subsec}, but which is trivial in
light of the previous theorem.  Namely, applying the theorem to the
prime pairs $(1,3)$ and $(4,6)$, we see that the deletion of the
dashed edges leaves a bipartite graph (with bipartition given by the
shading).  Theorem \ref{minimalprimepairs} now allows us to conclude
the original graph is invertible.
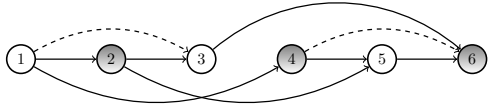
\begin{figure}[!ht]\label{dprimeex}
\begin{center}
\scalebox{.6}{
\begin{tikzpicture}
\GraphInit[vstyle=Normal]
\SetVertexMath
\SetVertexNoLabel

\Vertex[x=0.0cm,y=0.0cm]{a0}
\Vertex[x=2.0cm,y=0.0cm,style={shade}]{a1}
\Vertex[x=4.0cm,y=0.0cm]{a2}
\Vertex[x=6.0cm,y=0.0cm,style={shade}]{a3}
\Vertex[x=8.0cm,y=0.0cm]{a4}
\Vertex[x=10.0cm,y=0.0cm,style={shade}]{a5}

\AssignVertexLabel{a}{10}{$1$,$2$,$3$,$4$,$5$,$6$}
\tikzset{EdgeStyle/.append style = {->}}
\Edge(a0)(a1)
\Edge(a1)(a2)
\Edge(a3)(a4)
\Edge(a4)(a5)
\Edge[style={bend left=40}](a2)(a5)
\tikzset{EdgeStyle/.style = {->,bend right}}
\Edge(a0)(a3)
\Edge(a1)(a4)

\tikzset{EdgeStyle/.style = {->,dashed,bend left}}
\Edge(a0)(a2)
\Edge(a3)(a5)

\end{tikzpicture}}
\end{center}
\end{figure}
\end{Exam}

\noindent In the other direction, we wish to improve the ``lower
bound'' $\Gamma_D\leq H_D$, i.e., to find a supergraph of $\Gamma_D$
whose bipartiteness is forced by the invertibility of $D$.

\begin{Theo}\label{DeltaD}
Let $\Delta_D$ be the subgraph of $D$ obtained by removing every edge
$e_{ij}$ connecting vertices $i$ and $j$ with $\<i,j\>=1$.  Then
$\Delta_D$ is bipartite if $D$ is invertible.
\end{Theo}
\begin{proof}
If $D$ is invertible, then $\Gamma_D$ is bipartite.  Let $\l(i,j)$
denote the length of the maximal path from $i$ to $j$ in $D$ (or
$\Gamma_D$), and note that $\Gamma_D$ is formed from $D$ by removing
any edge $e_{ij}$ not on a maximal-length path from $i$ to $j$, i.e.,
removing $e_{ij}$ if $\l(i,j)>1$.  Since $\Delta_D$ is formed from $D$
by removing edges $e_{ij}$ with $\l(i,j)$ even, one can alternatively
view $\Delta_D$ as being constructed as the supergraph of $\Gamma_D$
where one adds back in an edge $e_{ij}$ from $D$ if and only if
$\l(i,j)$ is odd and greater than 1.  But adding an edge between
vertices with $\<i,j\>=-1$ can never construct an odd cycle, so the
bipartiteness of $\Gamma_D$ forces the bipartiteness of $\Delta_D$ as
well.
\end{proof}
\begin{Remark}
Analogously to Remark \ref{proper}, we note that if $D$ is invertible
and not itself bipartite, then $\Delta_D$ is a \emph{proper} bipartite
supergraph of $\Gamma_D$.
\end{Remark}

As above, it is of clear interest to ask if these bounds can be
tightened so as to uniquely identify a canonical subgraph of a general
$D$ whose bipartiteness is equivalent to the invertibility of $D$.
Without attempting to make this question particularly precise, we note
that the proof of the previous theorem convincingly prohibits such a
subgraph from existing.  Namely, $\Delta_D$ is quite clearly maximal
with respect to the property given in the theorem -- any supergraph of
$\Delta_D$ in $D$ contains an edge $e_{ij}$ between a pair of vertices
with $\<i,j\>=1$ and hence is not bipartite.  Thus the only candidates
for an ``optimal bipartite subgraph'' are subgraphs of $\Delta_D$, so
a single example of a non-invertible graph with bipartite $\Delta_D$
proves the impossibility of a subgraph whose bipartiteness is
equivalent to the invertibility of $D$.  We furnish such an example
below:

\begin{Exam}
Let $D$ be the graph in Figure \ref{DeltaDex}.  Then $D$ is
non-invertible, since the prime pair $(1,4)$ is not signable (we have
$\<1,4\>B_{1,4}^{-1}=(-1)(2-1)=-1$), but $\Delta_D$, obtained by
removing the dashed edges, is clearly bipartite.  We also note that
$\Delta_D$ is a proper supergraph of $\Gamma_D$ in this example, since
the edge $e_{2,5}\in \Delta_D-\Gamma_D$.
\begin{figure}[!ht]
\begin{center}
\scalebox{.6}{
\begin{tikzpicture}
\GraphInit[vstyle=Normal]
\SetVertexMath
\SetVertexNoLabel

\Vertex[x=0.0cm,y=0.0cm]{a0}
\Vertex[x=2.0cm,y=0.0cm,style={shade}]{a1}
\Vertex[x=4.0cm,y=0.0cm]{a2}
\Vertex[x=6.0cm,y=0.0cm,style={shade}]{a3}
\Vertex[x=8.0cm,y=0.0cm]{a4}

\AssignVertexLabel{a}{10}{$1$,$2$,$3$,$4$,$5$}
\tikzset{EdgeStyle/.append style = {->}}
\Edge(a0)(a1)
\Edge(a1)(a2)
\Edge(a2)(a3)
\Edge(a3)(a4)
\Edge[style={bend right=40}](a1)(a4)
\tikzset{EdgeStyle/.style = {->,dashed,bend left}}
\Edge(a0)(a2)
\Edge(a1)(a3)

\end{tikzpicture}}
\end{center}
\caption{A non-invertible graph with bipartite $\Delta_D$}\label{DeltaDex}
\end{figure}
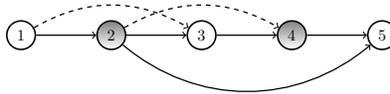
\end{Exam}

%%%%%%%%%%%%%%%%%%%%%%%%%%%%%%%%%%%%%%%%%%%%%%%%%%%%%%%%%%%%%%%%%%%%

\section{Constructing invertible bipartite graphs with unique perfect matching}

We begin by introducing some notation.  For vertices $d,d'\in D$,
write $d\ra d'$ if there is a directed edge from $d$ to $d'$, and
$d\lt d'$ if there is a directed path from $d$ to $d'.$ Recall that by
our conventions on labeling, our directed graphs have the property
that if $d\ra d'$, then $d<d'.$ In this section we present an
algorithm for constructing all invertible digraphs, and hence all
invertible bipartite graphs with unique perfect matching.  We begin
with the following straight-forward proposition:

\begin{Prop}\label{Construct}
All digraphs having the property that $d<d'$ whenever $d\ra d'$ can be
constructed from the graph on a single vertex (labelled $1$), via a
sequence of operations of the form:
\begin{quote}
If $D$ has vertices $\{1,2,\ldots,n-1\}$, choose a subset
$S_n\subseteq\{1,2,\ldots,n-1\}$ and add a vertex labelled $n$ and
(possibly multiple) adjacencies $s\ra n$ for all $s\in S_n$.
\end{quote}
\end{Prop}

To construct all \emph{invertible} digraphs, we find necessary and
sufficient conditions on the set $S$ in the proposition to ensure
invertibility at each stage of the construction process.  This is
clearly a requirement for the invertibility of the resulting graph.
If any intermediate graph $D$ of the construction process is
non-invertible then either the associated graph $\Gamma_D$ is not
bipartite or $D$ contains an unsignable prime pair.  Since neither of
these deficiencies can be rectified by the above process of adding
further vertices and/or edges, if at any stage of the construction
process the resulting graph is non-invertible, the final resulting
graph will be non-invertible.

\begin{Defn}
Given a subset $S$ of (the vertices of) an invertible directed graph
$D$ (on vertex set $\{1,2,\ldots,n-1\}$), let $D_S$ denote the graph
obtained by adjoining a new vertex labeled $n$ to $D$ and adding edges
$s\ra n$ for each $s\in S.$ Call $t\in S$ \emph{terminal} if $t\not\lt
s$ for all $s>t$ in $S$.  Call $S$ \emph{valid} if the maximal-path
subgraph $\Gamma_{D_S}$ is bipartite.
\end{Defn}

\begin{Lemma}
A subset $S$ of an invertible digraph $D$ is valid if and only if
$\langle t,n\rangle=-1$ in $D_S$ for all terminal $t\in S$.
\end{Lemma}

\begin{proof}
Since $D$ is invertible, $\Gamma_D$ is bipartite.  If $t$ is terminal,
then there exists a unique path from $t$ to $n$ in $D_S$, namely the
length-one path $t\ra n$.  Since this path necessarily survives to the
maximal-path subgraph, we must have $\langle t,n\rangle=-1$.  On the other hand, if $s$
is non-terminal, then the edge $s\ra n$ is not an edge of the maximal-path subgraph, and thus has no bearing on the bipartiteness of
$\Gamma_{D_S}$.
\end{proof}

\begin{Theo}\label{IneqThm}
Let $S$ be a valid subset of an invertible graph $D$.  For a vertex
$d\in D$, define $S_d=\{s\in S\,|\, d\leadsto s\}$, and partition $S_d$
into two sets $S_d^+$ and $S_d^-$ according to whether $\langle s,n\rangle=\pm
1$.  Then $D_S$ is invertible iff for all $d\in D$ with $|S_d|\geq 2$,
we have
\begin{align}\label{Ineq}
\sum_{s\in S_d^-}|B_{d,s}^{-1}|\geq \sum_{s\in S_d^+}|B_{d,s}^{-1}|.
\end{align}
\end{Theo}
\begin{proof}
Since $D$ is invertible, all prime pairs not including $n$ are already
signable, so by Theorem \ref{Alg}, it suffices to check the prime pairs
involving $n$, namely those prime pairs $(d,n)$ with $|S_d|\geq 2$.  Write
$S_d=\{s_1,\ldots,s_k\}$.  Repeatedly using the fact that $s_{j-1}$ is terminal in $D-\{s_j,s_{j+1},\ldots,s_k\}$, we have

\begin{align*}
B_{d,n}^{-1}&=B_{d,s_k}^{-1}B_{s_k,n}^{-1}+{_{s_k}}B_{d,n}^{-1}\\
&=-B_{d,s_k}^{-1}+{_{s_k}}B_{d,s_{k-1}}^{-1}\cdot {_{s_k}}B_{s_{k-1},n}^{-1}+{_{s_k,s_{k-1}}}B_{d,n}^{-1}\\
&=\qquad\vdots\\
&=-B_{d,s_k}^{-1}-B_{d,s_{k-1}}^{-1}-\cdots-B_{d,s_1}^{-1}\\
&=-\sum_{j=1}^k B_{d,s_j}^{-1}.
\end{align*}
Now use
$\langle d,n\rangle B_{d,s_j}^{-1}=\langle d,s_j\rangle\langle s_j,n\rangle B_{d,s_j}^{-1}=\langle s_j,n\rangle|B_{d,s_j}|$
by signability of $(d,s_j)$ to get
\begin{align*}
\langle d,n\rangle B_{d,n}^{-1}=-\langle d,n\rangle\sum_{s\in S_d}B_{d,s}^{-1}=-\sum_{s\in
  S_d}\langle d,s\rangle\langle s,n\rangle B_{d,s}^{-1}=-\sum_{s\in
  S_d}\langle s,n\rangle|B_{d,s}^{-1}|=\sum_{s\in S_d^-}|B_{d,s}^{-1}|-\sum_{s\in
  S_d^+}|B_{d,s}^{-1}|,
\end{align*}
whose positivity is the statement of the theorem.
\end{proof}

\begin{Cor}\label{monochrom}
Let $D$ and $S$ be as in the theorem.  If $D$ can be colored so that
$S$ is monochromatic, then $D_S$ is invertible.
\end{Cor}
\begin{proof}
We have $\langle s,n\rangle=-1$ for all $s\in S$, so $S_d^+$ is empty and the
inequality \eqref{Ineq} is trivially satisfied.
\end{proof}

\noindent As a corollary, we can immediately recover again the result of Godsil:
\begin{Cor}
If $D$ is bipartite, it is invertible.
\end{Cor}
\begin{proof}
If $D$ is bipartite, it can be constructed by a sequence of iterations
of the construction given above, adding a new right-most vertex and
connecting to a subset $S$ of the vertices, all of the same color.
\end{proof}

\begin{Cor}
Let $D$ and $S$ be as in the theorem.  If $|S|=1$, then $D_S$ is
invertible.
\end{Cor}
\begin{proof}
Since $|S|=1$, there are no $d\in D$ with $|S_d|\geq 2$, so the
condition is vacuously satisfied.
\end{proof}

\noindent Non-monochromatic sets $S$ do not appear to admit
particularly simple conditions for guaranteeing the invertibility of
$D_S$.  The next proposition deals with the case $|S|=2$.

\begin{Prop}\label{2setprop}
Let $S=\{s,s'\}$ be a valid subset of an invertible graph $D$.  Then
$D_S$ is invertible unless (and only unless) $\langle s,s'\rangle=-1$, $s\leadsto
s'$, and there exists $d\in D$ path-connected to both $s$ and $s'$
such that
$$
|B_{d,s'}^{-1}|< |B_{d,s}^{-1}|.
$$
In particular, if $s=1$, then $D_S$ is invertible.
\end{Prop}
\begin{proof}
First, if $\langle s,s'\rangle=1$, then $D_S$ is invertible by
Corollary \ref{monochrom}.  We thus assume that $\langle
s,s\rangle=-1$.  Further, we may assume that $s\leadsto s'$, since
otherwise the subset $\{s,s',n\}\subset \Gamma_{D_S}$ could not be
2-colored, contradicting our assumption that $S$ was valid.  Since $D$
is invertible, it suffices to check the new prime pairs created upon
addition of $n$, which are exactly the prime pairs $(d,n)$ for $d$
such that $d\leadsto s$.  For such $d$, we have $S_d=\{s,s'\}$,
independently of $d$.  Since $\langle s,n\rangle=1$, and $\langle
s',n\rangle=-1$, Theorem \ref{IneqThm} gives that $D_S$ is invertible
if and only if $|B_{d,s'}^{-1}|\geq |B_{d,s}^{-1}|$ for each such $d$.
\end{proof}

\begin{Cor}
Let $S=\{s,s'\}$ be a valid subset of an invertible graph $D$.  Then
there exists a positive integer $k$ such that if $k$ edges are added
to $D$ connecting $s$ to $s'$, then $D_S$ is invertible.
\end{Cor}
\begin{proof}
From the proposition, $D_S$ is invertible if for all $d$
path-connected to both $s$ and $s'$, we have $|B_{d,s'}^{-1}|\geq
|B_{d,s}^{-1}|$.  Writing
$B_{d,s'}^{-1}=B_{d,s}^{-1}B_{s,s'}^{-1}+{_s}B_{d,s'}^{-1}$ and
canceling a $B_{d,s}^{-1}$ from both sides (if $B_{d,s}^{-1}=0$, the
inequality is trivially satisfied), we arrive
at the invertibility condition
$$
\left|\frac{{_s}B_{d,s'}^{-1}}{B_{d,s}^{-1}}+B_{s,s'}^{-1}\right|\geq 1.
$$ Since $B_{s,s'}^{-1}$ is independent of $d$, one can ensure that this
inequality is satisfied for all such $d$ by choosing, for example,
$$ B_{s,s'}^{-1}\geq 1+\max_{d:
  B_{d,s}^{-1}\neq0}\left|\frac{{_s}B_{d,s'}^{-1}}{B_{d,s}^{-1}}\right|.
$$
\end{proof}

%%%%%%%%%%%%%%%%%%%%%%%%%%%%%%%%%%%%%%%%%%%%%%%%%%%%%%%%%%

\subsection{Application:  Constructing invertible unicyclic graphs}\label{unicycleconstruct}

As an extended example of the iterative construction process, we turn
our attention to directed graphs $D$ which are associated to bipartite
unicyclic graphs with a unique perfect matching (and in particular,
invertible such digraphs).  As in the previous section, we begin with
a single vertex $D_1$ and iteratively adjoin edges to new vertices.
More precisely, we proceed inductively for $i\geq 2$, letting $D_i$ be
the graph obtained by adding a new vertex $i$ to $D_{i-1},$ and edges
$s\ra i$ for all $s$ in some \emph{adjacency set} $S_i\subseteq
[1,i-1].$ Given adjacency sets $S_1,\ldots,S_{n}$, we will often
denote simply by $D$ rather than $D_n$ the end result of iteratively
adjoining each $S_i$.

\medskip We break the construction of all unicyclic graphs into two
steps: first the construction of the cycle, and then the rest of the
graph.  The idea in both steps is to precisely describe the
combinatorial restrictions on the number of edges added at each stage
of the iteration (i.e., on $|S_i|$) forced by unicyclicity.  To this
end, we introduce the following terminology:

\begin{Defn}
A \emph{Motzkin partition} of a positive integer $N$ is a partition
$P=\{p_1,\ldots,p_N\}$ of $N$ into exactly $N$ parts such that $0\leq
p_i\leq 2$ for all $i$, and $\sum\limits_{k=1}^i p_k<i$ for all $i<N$.
Note that the partial sum condition forces $p_1=0$ and $p_N=2$.
\end{Defn}
\begin{Remark}
We so-name these partitions due to their connection to the well-known
Motzkin numbers (e.g., \cite{OEIS}, sequence A001006).  The authors
wish to thank David Speyer \cite{MO10} for pointing out this link to
us.  We note that a previous combinatorial interpretation of Motzkin
numbers in terms of unicyclic graphs does not seem to exist in the
literature (see, e.g., \cite{DS77}).
\end{Remark}

Continuing our slight abuse of terminology from Section \ref{UInvert},
we say that a digraph $D$ is the cycle graph (resp. unicyclic or
acyclic) if it is the cycle graph (resp. unicyclic or acyclic) as an
undirected graph.  We begin by describing the possible iterative
constructions of the cycle graph on $N$ vertices.  To avoid
conflicting with future notation, we use $C_i$ instead of $S_i$ to
describe the adjacency sets for the cycle.

\begin{Prop}\label{uconst}
Let $P$ be a Motzkin partition of $N$.  Then there exist adjacency
sets $C_i\subseteq[1,i-1]$ for $2\leq i\leq N$, such that $|C_i|=p_i$
and such that the digraph $D$ resulting from iteratively adjoining
these adjacency sets is the cycle graph on $N$ vertices.
\end{Prop}

\begin{proof}
We will construct sets $C_i$ for $2\leq i\leq n$ such that $|C_i|=p_i$
and such that the resulting graph $D$ is two-regular and connected.
From this it easily follows that, as an undirected graph, $D$ is the
cycle on $\sum p_i=N$ vertices.  As always, let $D_1$ denote the
single-vertex graph.  We proceed inductively, for $i\geq 2$ choosing
$C_i\subset[1,i-1]$ of size $p_i$ and letting $D_i$ be obtained from
$D_{i-1}$ by adding vertex $i$ and edges $c\to i$ for each $c\in C_i$.
Note that to prove $D$ is two-regular and connected it is necessary
and sufficient to prove that we can chose our adjacency sets such that
$D_i$ is acyclic for all $i<N,$ and has maximum vertex degree at most
$2$ for all $i.$ These properties are clearly satisfied by $D_1.$ Now
the induction: Suppose that $D_i$ ($1\leq i\leq N-2$) is acyclic with
maximum vertex degree at most 2.  We show that we can choose $C_{i+1}$
so that $|C_{i+1}|=p_{i+1}$ and $D_{i+1}$ again is acyclic with
maximum vertex degree at most 2.  If $p_{i+1}=0$, this is trivial.  If
$p_{i+1}=1$, we simply need to prove the existence of a vertex of
degree at most 1 in $D_i$: If such a vertex did not exist, we would
have $\sum_{k=1}^i p_k=i$, contradicting the partial sum condition of
Motzkin partitions.  Similarly, if $p_{i+1}=2$ and there are not two
unconnected vertices of degree at most one, then $D_i$ is the path
graph on $i$ vertices and we have $\sum_{k=1}^i p_k=i-1$.  But now
$\sum_{k=1}^{i+1}p_k=i-1+2=i+1$, again violating the partial sum
condition.  Thus for all $i<N$ we can select the desired adjacency
sets $C_i$.
\end{proof}

Next we show how to iteratively construct all graphs with a prescribed
number of vertices $M$ and a unique cycle on a prescribed set of $N$
vertices.  From this it is easy to see that we can so construct all
unicyclic graphs (Remark \ref{allunis}).  Recall that by choosing a
different topological sort if necessary, we can assume that the cycle
occurs on consecutive vertices, say on vertices $\{v+1,\ldots,v+N\}$.
The apparent combinatorial complexity of the following proposition is
a result of having to chose adjacency sets $S_i$ so that they include
the sets $C_i$ as above (thus constructing a cycle of the desired
length), while ensuring that a second cycle is never created.

\begin{Prop}\label{uconst2}
Let $M$ and $N$ be positive integers with $M>N$, and let $v$ be an
integer between $1$ and $M-N$.  Let $P=\{p_{v+1},\ldots,p_{v+N}\}$ be
a Motzkin partition of $N$ and $T=\{t_1,\ldots,t_M\}$ be a partition
of $M$ into $M$ parts such that
\begin{itemize}
\item $t_i\geq p_i$ for all $v+1\leq i\leq v+N;$
\item $\sum_{k=1}^i t_k\leq i$ for all $i$;
\item $\sum_{k=1}^i t_k<i$ for all $1\leq i< v+N$.
\end{itemize}
Then there exist adjacency sets $S_i\subset [1,i-1]$ of size $t_i$,
for $1< i\leq M$, such that the digraph $D$ resulting from the
iterative construction process has $M$ vertices, $M$ edges, and a
single undirected cycle of length $N$ on vertices $v+1$ through $v+N.$
\end{Prop}

\begin{proof}
For the proof of this proposition, we consider a slight alteration to
the construction process.  Instead of adjoining a new vertex $i$ at
the $i$-th step, we begin with the completely disconnected graph on
$M$ vertices and, for $1<i\leq M$, consider adding the directed edges
$s\ra i$ for all $s$ in the adjacency set $S_i$.  Define the partition
$Q=\{q_1,\ldots,q_M\}$ of $M-N$ into $M$ parts by
\[
q_i=
\begin{cases}
t_i-p_i&\text{ for }v+1\leq i\leq v+N\\
t_i&\text{ else}
\end{cases}.\]
By Proposition \ref{uconst}, for $v+1<i\leq v+N$ there exist adjacency
sets $C_i$ such that $|C_i|=p_i$ and the graph resulting from adding
the edges $c\ra i$ for all $c$ in $C_i$ and all $v+1<i\leq v+N$ is the
cycle graph on vertices $\{v+1,\ldots,v+N\}.$ We will prove the
existence of sets $S_i$ ($1<i\leq M$), by showing the existence of
sets $R_i$ with $|R_i|=q_i$, and setting $S_i=R_i\cup C_i$ (here
$C_i=\emptyset$ for all $i\not\in [v+2,v+N]$).  Create the graph on
$M$ vertices consisting of a single cycle on vertices
$\{v+1,\ldots,v+N\}$ by adding the directed edges $c\ra i$ prescribed
by the adjacency sets $C_i$ for $v+1<i\leq v+N.$ Call this graph
$H_1.$ For $1<i\leq M,$ let $H_i$ denote the graph obtained from
$H_{i-1}$ by adding edges $s\ra i$ for all $s\in R_i.$ Assuming the
graph $H_{i-1}$ is unicyclic, it remains to show that in the $i$-th
step we can select the $q_i$ vertices of $R_i$ such that the resulting
graph is unicyclic., i.e., that the part of the graph to the left of
$i$ has at least $q_i$ connected components.  For ease of notation, we
split into several cases.

\medskip

\noindent \emph{Case I $(1< i\leq v)$:}  Let $h_i$ be the number of
connected components of $H_{i-1}\cap[1,i-1]$.  Since there are no
cycles in $H_{i-1}\cap[1,i-1]$, we have $h_i=i-1-\sum_{k=1}^{i-1} q_k$.
The condition that there are sufficiently many connected components in
constructing $H_i$ is simply $q_i\leq h_i$, i.e., that
$\sum_{k=1}^iq_k\leq i-1<i$, but this is true by assumption (noting $q_k=t_k$
for $k<v$).

\medskip

\noindent \emph{Case II $(v+1\leq i\leq v+N)$:} Since any additional
edge between two vertices of the cycle would force a second cycle, for
vertices in the range $[v+1,v+N]$, it suffices to show that for each
such $i$, the number of connected components in $H_{i-1}\cap[1,v+N]$,
given by $h_i=v+1-\sum_{k=1}^{i-1} q_k,$ is at least $q_i+1$ (the
cycle accounts for the extra connected component).  Again, this is
equivalent to the condition that $\sum_{k=1}^i q_k\leq v$, which
follows from

\[
\sum_{k=1}^{i} q_k\le\sum_{k=1}^{v+N}q_k=\sum_{k=1}^{v+N}t_k-N\le v+N-N.
\]

\medskip

\noindent \emph{Case III $(v+N+1\leq i\leq M)$:} For $i>v+N$, the
number of connected components in $H_{i-1}\cap[1,i-1]$ is calculated
by considering the $v$ initial vertices, the cycle, and the
$i-(v+N)-1$ vertices between $v+N$ and $i$ (i.e. the number of
vertices between the cycle and $i$).  Each of the
$\sum_{k=1}^{i-1}q_k$ added edges joins two of these components.  Thus
the number of connected components in $H_{i-1}\cap[1,i-1]$ is
\[
h_i=v+1+i-(v+N)-1-\sum_{k=1}^{i-1}q_k=i-N-\sum_{k=1}^{i-1}q_k,
\]
and the unicyclicity condition is just $q_i\leq
i-N-\sum_{k=1}^{i-1}q_k.$ Equivalently, we need $\sum_{k=1}^i q_k\le
i-N$, which follows from
\[
\sum_{k=1}^i q_k=\sum_{k=1}^i t_k-N\leq i-N.
\]
\end{proof}

\noindent We illustrate this construction process with the following example.

\begin{Exam}
Let $M=8$ and $N=5$ and $v=2.$ The following partitions meet the
requirements of Proposition \ref{uconst} and
\ref{uconst2}: $$T=\{0,1,0,0,2,2,2,1\}\qquad P=\{0,0,2,1,2\}.$$ To
begin, we construct the cycle on vertices $\{3,4,5,6,7\}$.  By the
requirement that $|C_i|=p_i$, the fact that $p_3=p_4=0$ forces the
adjacency sets $C_3$ and $C_4$ to be empty.  The only possible choices
for $C_5$ is then $\{3,4\}$.  There are two valid options for the set
$C_6,$ namely $\{3\}$ or $\{4\}$.  $C_7$ is determined by our choice
for $C_6,$ with $C_7=\{4,6\}$ if $C_6=\{3\}$ and $C_7=\{3,6\}$
otherwise.  Thus, there are two possibilities for the intermediate
graph $H_1$ consisting of just the cycle, as shown in Figure
\ref{constexam0}.

\begin{center}
\begin{figure}[!ht]\label{constexam0}
\begin{tabular}{ccc}
\scalebox{.4}{
\begin{tikzpicture}
\GraphInit[vstyle=Normal]
\SetVertexMath
\SetVertexNoLabel

\Vertex[x=0.0cm,y=0.0cm]{a0}
\Vertex[x=2.0cm,y=0.0cm]{a1}
\Vertex[x=4.0cm,y=0.0cm]{a2}
\Vertex[x=6.0cm,y=0.0cm]{a3}
\Vertex[x=8.0cm,y=0.0cm]{a4}
\Vertex[x=10.0cm,y=0.0cm]{a5}
\Vertex[x=12.0cm,y=0.0cm]{a6}
\Vertex[x=14.0cm,y=0.0cm]{a7}

\AssignVertexLabel{a}{10}{$1$,$2$,$3$,$4$,$5$,$6$,$7$,$8$}
\tikzset{EdgeStyle/.append style = {->}}
\Edge(a3)(a4)
\Edge(a5)(a6)
\tikzset{EdgeStyle/.style = {->,bend left}}
\Edge(a2)(a4)
\renewcommand*{\EdgeLineWidth}{1mm}
\tikzset{EdgeStyle/.style   = {->,bend left = 40, line width = \EdgeLineWidth}}
\Edge(a2)(a5)
\tikzset{EdgeStyle/.style = {->,bend right, line width = \EdgeLineWidth}}
\Edge(a3)(a6)
\end{tikzpicture}}
&
\hspace*{.4in}
&
\scalebox{.4}{
\begin{tikzpicture}
\GraphInit[vstyle=Normal]
\SetVertexMath
\SetVertexNoLabel

\Vertex[x=16.0cm,y=0.0cm]{a0}
\Vertex[x=18.0cm,y=0.0cm]{a1}
\Vertex[x=20.0cm,y=0.0cm]{a2}
\Vertex[x=22.0cm,y=0.0cm]{a3}
\Vertex[x=24.0cm,y=0.0cm]{a4}
\Vertex[x=26.0cm,y=0.0cm]{a5}
\Vertex[x=28.0cm,y=0.0cm]{a6}
\Vertex[x=30.0cm,y=0.0cm]{a7}

\AssignVertexLabel{a}{10}{$1$,$2$,$3$,$4$,$5$,$6$,$7$,$8$}
\tikzset{EdgeStyle/.append style = {->}}
\Edge(a3)(a4)
\Edge(a5)(a6)
\tikzset{EdgeStyle/.style = {->,bend left}}
\Edge(a2)(a4)
\renewcommand*{\EdgeLineWidth}{1mm}
\tikzset{EdgeStyle/.style   = {->,bend right, line width = \EdgeLineWidth}}
\Edge(a3)(a5)
\Edge[style={bend right=40}](a2)(a6)
\end{tikzpicture}}
\end{tabular}
\caption{Two possible intermediate digraphs $H_1$ resulting from the
  iterative construction process on partition $P=\{0,0,2,1,2\},$ with
  $M=8$ and $v=2.$}
\end{figure}
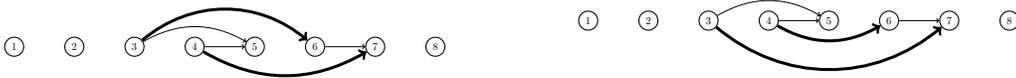
\end{center}

Since $|S_i|=t_i$, $S_3,$ and $S_4$ are the empty set.  The
only possible non-empty adjacency set for $S_2$ is $\{1\}.$ For $3\le
i\le 7,$ $C_i\subseteq S_i$.  By our definitions of $P$ and $T$,
$|S_5|=|C_5|, |S_7|=|C_7|,$ and $|S_6|=|C_6|+1.$ Thus, sets $S_5$ and
$S_6$ are completely determined by $H_1$, while $S_6$ contains the
subset $C_6$ and one additional vertex chosen from $\{1,2\}$.  The
final adjacency set $S_8$ can equal any single vertex in the set
$\{1,2,\ldots,7\}$.  Figure \ref{constexam2} shows the two possible
resulting digraphs if we set $S_6=C_6\cup\{2\}$ and $S_8=\{7\}$.

\begin{center}
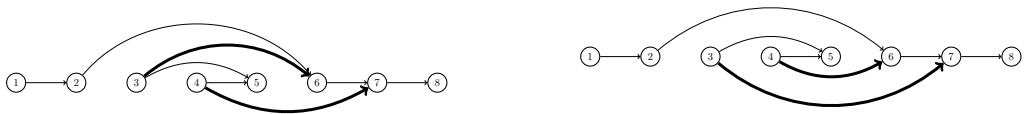
\begin{figure}[!ht]
\begin{tabular}{ccc}
\scalebox{.4}{
\begin{tikzpicture}
\GraphInit[vstyle=Normal]
\SetVertexMath
\SetVertexNoLabel

\Vertex[x=0.0cm,y=0.0cm]{a0}
\Vertex[x=2.0cm,y=0.0cm]{a1}
\Vertex[x=4.0cm,y=0.0cm]{a2}
\Vertex[x=6.0cm,y=0.0cm]{a3}
\Vertex[x=8.0cm,y=0.0cm]{a4}
\Vertex[x=10.0cm,y=0.0cm]{a5}
\Vertex[x=12.0cm,y=0.0cm]{a6}
\Vertex[x=14.0cm,y=0.0cm]{a7}

\AssignVertexLabel{a}{10}{$1$,$2$,$3$,$4$,$5$,$6$,$7$,$8$}
\tikzset{EdgeStyle/.append style = {->}}
\Edge(a0)(a1)
\Edge(a3)(a4)
\Edge(a5)(a6)
\Edge(a6)(a7)
\tikzset{EdgeStyle/.style = {->,bend left}}
\Edge(a2)(a4)
\Edge[style={bend left=50}](a1)(a5)
\renewcommand*{\EdgeLineWidth}{1mm}
\tikzset{EdgeStyle/.style   = {->,bend left = 40, line width = \EdgeLineWidth}}
\Edge(a2)(a5)
\tikzset{EdgeStyle/.style = {->,bend right, line width = \EdgeLineWidth}}
\Edge(a3)(a6)
\end{tikzpicture}}
&
\hspace*{.4in}
&
\scalebox{.4}{
\begin{tikzpicture}
\GraphInit[vstyle=Normal]
\SetVertexMath
\SetVertexNoLabel

\Vertex[x=16.0cm,y=0.0cm]{a0}
\Vertex[x=18.0cm,y=0.0cm]{a1}
\Vertex[x=20.0cm,y=0.0cm]{a2}
\Vertex[x=22.0cm,y=0.0cm]{a3}
\Vertex[x=24.0cm,y=0.0cm]{a4}
\Vertex[x=26.0cm,y=0.0cm]{a5}
\Vertex[x=28.0cm,y=0.0cm]{a6}
\Vertex[x=30.0cm,y=0.0cm]{a7}

\AssignVertexLabel{a}{10}{$1$,$2$,$3$,$4$,$5$,$6$,$7$,$8$}
\tikzset{EdgeStyle/.append style = {->}}
\Edge(a0)(a1)
\Edge(a3)(a4)
\Edge(a5)(a6)
\Edge(a6)(a7)
\tikzset{EdgeStyle/.style = {->,bend left}}
\Edge(a2)(a4)
\Edge[style={bend left=40}](a1)(a5)
\renewcommand*{\EdgeLineWidth}{1mm}
\tikzset{EdgeStyle/.style   = {->,bend right, line width = \EdgeLineWidth}}
\Edge(a3)(a5)
\Edge[style={bend right=40}](a2)(a6)
\end{tikzpicture}}
\end{tabular}
\caption{Two possible unicyclic digraphs resulting from the iterative
  construction process on partitions $T=\{0,1,0,0,2,2,2,1\}$ and
  $P=\{0,0,2,1,2\}.$}\label{constexam2}
\end{figure}
\end{center}

\end{Exam}

\begin{Remark}\label{allunis}
Every $D$ which is connected and unicyclic as an undirected graph can
be constructed in this manner.  Let $D$ be such a graph, with $M$
vertices and a cycle of length $N$.  As above, without loss of
generality we assume that the vertices contained in the unique cycle
are labelled consecutively, say $v+1$ through $v+N$.  For each vertex
$i$ in $D$, let $S_i\subset[1,i-1]$ be the set of vertices
edge-connected to $i$ and for $i\in[v+1,v+N]$, let $C_i\subset S_i$
correspond to those edges coming from other vertices in the cycle.
Then by connectedness and unicyclicity, the sets $T=\{|S_i|\}$ and
$P=\{|C_i|\}$ form partitions as in the proposition, and we can follow
the construction process in the proof of Proposition \ref{uconst2}
with adjacency sets $S_i$ and $C_i$ to reconstruct $D$.
\end{Remark}

\begin{Cor}
Given the construction process described in the proof of Proposition
\ref{uconst2}, the resulting digraph $D$ is invertible iff $N$ is even
or $C_i=\{i-1\}$ for all $v+1<i<v+N$.
\end{Cor}
\begin{proof}
If $N$ is even, the resulting unicyclic digraph is bipartite, and thus
invertible.  The requirement that $C_i=\{i-1\}$ for all $v+1<i<v+N$
forces $C_{v+N}=\{v+1, v+N-1\},$ and is equivalent to saying the
undirected cycle subgraph $D'$ induced by vertices
$\{v+1,\ldots,v+N\}$ must contain only one source and one sink vertex
(i.e., $k=1$ in the language of Theorem \ref{unicycleinvert}), and
that these two vertices must be adjacent.  Thus, the result follows
directly from Theorem \ref{unicycleinvert}.
\end{proof}

%%%%%%%%%%%%%%%%%%%%%%%%%%%%%%%%%%%%%%%%%%%%%%%%%%%%%%%%%%%%%%%%%%%%%%%%%%%%%%%%%%%%%%%%%%

\end{document}